%% file: filter-02-11-nofourier-arXiv-v2.tex
\documentclass[11pt,a4paper]{article}
\usepackage{fullpage}
\usepackage{microtype}
\usepackage{graphicx}

\input{preamble}

\newcommand{\hspn}{\hspace{-0.2cm}}
\def\diag{{\mathrm{diag}}}

\begin{document}

\title{\fontdimen2\font=4.5pt\text{Adaptive~Denoising~of~Signals~with~Local~Shift-Invariant~Structure}}
\author{
Zaid Harchaoui\thanks{University of Washington, Seattle, WA 98195, USA. Email: \texttt{zaid@uw.edu}.}
	\and Anatoli Juditsky\thanks{LJK, Universit\'e Grenoble Alpes, 700 Avenue Centrale,  38401 %Domaine Universitaire de
	Saint-Martin-d'H\`eres, France. Email: \texttt{anatoli.juditsky@univ-grenoble-alpes.fr}.}
	\and Arkadi Nemirovski\thanks{Georgia Institute of Technology, Atlanta, GA 30332, USA. Email: \texttt{nemirovs@isye.gatech.edu}.}
\and Dmitrii Ostrovskii\thanks{University of Southern California, Los Angeles, CA 90089, USA. Email: \texttt{dostrovs@usc.edu}.}
}

\maketitle

\begin{abstract}
We discuss the problem of adaptive discrete-time signal denoising in the situation where the signal to be recovered admits a ``linear oracle''---an unknown linear estimate that takes the form of convolution of observations with a time-invariant filter.
It was shown by Juditsky and Nemirovski (2009) \cite{jn1-2009} that when the $\ell_2$-norm of the oracle filter is small enough, such oracle can be ``mimicked'' by an efficiently computable \textit{adaptive} estimate of the same structure with the observation-driven filter.
The filter in question was obtained as a solution to the optimization problem in which the $\ell_\infty$-norm of the Discrete Fourier Transform (DFT) of the estimation residual is minimized under constraint on the $\ell_1$-norm of the filter DFT.
In this paper, we discuss a new family of adaptive estimates which rely upon minimizing the $\ell_2$-norm of the estimation residual.
We show that such estimators possess better statistical properties than those based on~$\ell_\infty$-fit;
in particular, under the assumption of \textit{approximate shift-invariance} we prove oracle inequalities for their $\ell_2$-loss and improved bounds for $\ell_2$- and pointwise losses.
We also study the relationship of the approximate shift-invariance assumption with the {\em signal simplicity} introduced in \cite{jn1-2009} and discuss the application of the proposed approach to {\em harmonic oscillation} denoising.
\end{abstract}

%\tableofcontents

\section{Introduction}
\label{sec:intro}
The problem we consider in this paper is that of signal denoising:
given noisy observations
\be
y_\tau = x_\tau + \sigma \zeta_\tau, \quad \tau \in \Z
\ee{eq:intro-observations}
we aim at recovering a signal $(x_t)_{t\in \Z}$.  It is convenient for us to assume that signal and noises are complex-valued.
Observation noises $\zeta_\tau$ are assumed to be independent of $x$ i.i.d.~standard complex-valued Gaussian random variables (denoted~$\zeta_{\tau} \sim \CN(0,1)$), meaning that $\zeta_\tau=\zeta^1_\tau+i\zeta^2_\tau$ with i.i.d.~$\zeta^1_\tau,\zeta^2_\tau\sim \N(0,1)$.
Our goal may be, for instance, to recover the value $x_t$ of the signal at time $t$ given observations $(y_\tau),\,{|\tau-t|\leq m}$ for some $m\in \Z_+$ (problem referred to as {\em signal interpolation} in signal processing literature), or to estimate the value $x_{t+h}$ given observations $(y_\tau),\,{t-m\leq \tau\leq t}$ ({\em signal prediction} or {\em extrapolation}), etc.
\par
The above problem is classical in statistics and signal processing. In particular, {\em linear estimates} of the form
\[
\wh{x}_t=\sum_{\tau\in \Z}\phi_\tau y_{t-\tau}
\]
are ubiquitous in nonparametric estimation; for instance, classical kernel estimators are of this type.
More generally, linear estimates are considered both theoretically attractive and easy to use in practice \cite{ibragimov1988,donoho1992,donoho1994,kailath2000linear,tsybakov_mono,wasserman2006all}.
When the set $\X$ of signals is well-specified, one can usually compute a (nearly) minimax on $\X$ linear estimator in closed form. In particular, if $\X$ is a class of ``smooth signals,'' such as a H\"older or a Sobolev ball, then the corresponding estimator is given by the kernel estimator with properly selected bandwidth~\cite{tsybakov_mono}, and is minimax among all possible estimators.
Moreover, linear estimators are known to be nearly minimax optimal with respect to the pointwise loss \cite{ibragimov1984,donoho1994} and the~$\ell_2$-loss~\cite{donoho1990minimax,pinsker1980,juditsky2018near-aos,juditsky2018near-msl} under rather general assumptions about the set $\X$ of possible signals.
Besides this, if the set $\X$ of signals is specified in a computationally tractable way, then a near-minimax linear estimator can be efficiently computed by solving a convex optimization problem~\cite{juditsky2018near-aos}, \cite{juditsky2018near-msl}.
\par
The strength of this approach, however, comes at a price:
%However, the above-described linear estimation approach has a crucial disadvantage:
in order to implement the estimate the set $\X$ must be {\em known to the statistician.} Such knowledge is crucial: near-minimax estimator for one signal set can be of poor quality for another one. Thus, linear estimation approach cannot be directly implemented when no prior knowledge of $\X$ is available.
In the statistical literature this difficulty is usually addressed via adaptive \textit{model selection}~\cite{birge1997model,goldenshluger2011bandwidth,johnstone-book,lepski1991,lepski1997optimal,lepski2015adaptive,massart2007concentration,tsybakov_mono}. However, model selection procedures usually impose strong structural assumptions on the signal set, assuming it to be known up to a few hyper-parameters.\footnote{More general adaptation schemes have been recently introduced, e.g.,  routines from \cite{goldenshluger2013general,lepski2015adaptive} which can handle, for example, adaptation to inhomogeneous and anisotropic smoothness of the signal. However, the proposed schemes cannot be implemented in a numerically efficient fashion, and therefore are not practical.}

An alternative approach to the denoising problem with unknown $\cX$ was proposed in~\cite{nemirovsky1991nonparametric}.
There, instead of directly restricting the class of signals and requiring a specification of $\cX$, one restricts the class of {possible estimators.}
Namely, let us denote $\C(\Z)$ the space of complex-valued functions on $\Z$, and let, for $m\in \Z_+$,~$\C_m(\Z)$ be the space of complex-valued sequences that vanish outside the set~$\{-m,...,m\}$. We consider linear {\em convolution-type} estimators, associated with {\em filters}~$\phi \in \C_m(\Z)$ of the form
\begin{equation}
\label{eq:convolution-estimate}
\wh{x}_t = [y * \phi]_t := \sum_{\tau \in \Z} \phi_\tau y_{t-\tau} = \sum_{|\tau|\leq m} \phi_\tau y_{t-\tau}.
\end{equation}
%Note that the summation in the above formula is in fact finite due to~$\phi \in \C_n(\Z)$.
%It was observed in \cite{jn1-2009} that
Informally, the problem we are interested in here is as follows:
\begin{quote}
{\em If we fix the structure~\rf{eq:convolution-estimate} of the estimate and consider {the form of the filter $\phi$} as a
``free parameter,'' is it possible to build an estimation procedure which is adaptive with respect to this parameter?}
\end{quote}
In other words,
suppose that a ``good'' filter $\phi^o$  with small estimation error ``exists in nature.''
Is it then possible to construct a data-driven estimation routine which has (almost) the same accuracy as the ``oracle''---a hypothetic optimal estimation method utilizing $\phi^o$?

The above question was first answered positively in~\cite{jn1-2009} using the estimation machinery from~\cite{nemirovsky1991nonparametric}.
To present the ideas underlying the approach of~\cite{jn1-2009} we need to define the class of ``well-filtered'' or ``simple'' signals~\cite{harchaoui2015adaptive,jn1-2009}.

\begin{definition}[{Simple signals}]
\label{def:l2simple}
Given parameters~$m, n \in \Z_+$,~$\rho \ge 1$, and~$\theta \ge 0$, signal $x \in \C(\Z)$ is called $(m, n, \rho,\theta)$-\emph{simple} if there exists $\phi^o \in \C_{m}(\Z)$ satisfying
\begin{equation}
\label{eq:l2-oracle-norm}
\|\phi^o\|_{2} \le {\rho \over \sqrt{2m+1}},
\end{equation}
and such that
\begin{equation}
\label{eq:l2-oracle-bias}
|x_\tau - [\phi^o*x]_{\tau}| \le {\sigma \theta\rho \over \sqrt{2m+1}}, \quad \text{for all} \;\;\; |\tau|\leq m+n.
\end{equation}
%We denote $\S_{m,n}(\rho,\theta)$ the set of $(m,n, \rho,\theta)$-simple signals for a given quadruple~$(m,n,\rho,\theta)$.
%The signals which are $(n, \rho,\theta)$-simple for given $\rho,\theta$, and any $n \ge \rho^2$, are called$(\rho,\theta)$-\emph{parametric}; the class of such signals is denoted $\S(\rho,\theta)$.
\end{definition}
%\noindent
Decomposing the pointwise mean-squared error of the estimate~\rf{eq:convolution-estimate} with $\phi=\phi^o$ as
\[
\E |x_\tau - [\phi^o * y]_\tau|^2 = \sigma^2 \E|[\phi^o * \zeta]_\tau|^2 + |x_\tau - [\phi^o * {x}]_\tau|^2,
\]
 we immediately arrive at the following bound on the pointwise expected error:
\begin{align}
\label{eq:ss1}
\left[\E |x_\tau - [\phi^o * y]_\tau|^2\right]^{1/2} \le {\sigma \sqrt{1+\theta^2}\rho\over \sqrt{2m+1}}, \quad |\tau|\leq {m+n}.
\end{align}
In other words, simple signals are those for which there exists a linear estimator (i) utilizing observations in the $m$-neighbourhood of a point, (ii) invariant in the $(m+n)$-vicinity of the origin, and (iii) attaining pointwise risk of order $m^{-1/2}$ in that vicinity (For brevity, here we refer to the quantity~$\E[|x_t - \wh x_t|^2]^{1/2}$ as the pointwise risk (at~$t \in \Z$) of estimate~$\wh x$. Parameters~$\rho,\theta$ allow for refined control of the risk and specify the bias-variance balance.

Now, assume that the only prior information about the signal to be recovered is that it %belongs to the class $\S_{m,n}(\rho,\theta)$
is $(m, n, \rho,\theta)$-simple with some known $(m,n,\rho,\theta)$.
As we have just seen, this implies existence of a convolution-type linear estimator~$\wh x^o = \phi^o * y$ with a good statistical performance.
The question is whether we can use this information to ``mimic'~$\wh x^o$, i.e.,~to construct an estimator of $(x_\tau)_{|\tau|\leq n}$ with comparable statistical performance when only using available observations.
Answering this question is not straightforward. To build an adaptive estimator, one could implement the cross-validation procedure by minimizing some observable proxy of the quadratic loss of the estimate, say, the $\ell_2$-norm of the residual $([y-\vphi * y]_{\tau} )_{|\tau|\leq m+n}$, over the set of filters $\varphi$ satisfying~\rf{eq:l2-oracle-norm}. However, it is well known that the set of filters satisfying~\rf{eq:l2-oracle-norm} is too ``massive'' to allow for construction of adaptive estimate with the risk bound similar to~\rf{eq:ss1} even when $\rho=1$.\footnote{While this statement appears self-evident to statisticians of older generations, younger researchers may expect an explanation. This is why we provide a brief discussion of the ``naive estimate'' in Section \ref{sec:naive} of the appendix.} As a result, all known to us approaches to adaptive estimation in this case impose some extra constraints on the class of filters such as regularity~\cite{efromovich1996sharp} or sparsity in a certain basis \cite{donoho1994ideal}, etc.

Nevertheless, surprisingly, adaptive convolution-type estimators with favorable statistical performance guarantees can be constructed.
The key idea, going back to~\cite{jn1-2009}, is to pass to a ``new oracle'' with a characterization which better suits the goal of adaptive estimation.
Namely, one can easily verify (cf., e.g.,~\cite[Proposition 3]{harchaoui2015adaptive}) that if a filter $\phi^o\in \C_m(\Z)$ satisfies relations~\rf{eq:l2-oracle-norm} and \rf{eq:l2-oracle-bias}, then its auto-convolution~$\vphi^o=\phi^o * \phi^o \in \C_{2m}(\Z)$ (with twice larger support) satisfies their analogues
\begin{eqnarray}
\label{eq:intro-autoconv-norm}
\left\|F_{2m} [\vphi^o] \right\|_1 \hspn &\le& \hspn \frac{2\rho^2}{\sqrt{4m+1}},\\
|x_\tau - [\vphi^o * x]_\tau| \hspn &\le& \hspn \frac{2\sqrt{2}\sigma\theta\rho^2}{\sqrt{4m+1}}, \quad |\tau|\leq n;\nonumber
\end{eqnarray}
here~$F_n$ is the unitary Discrete Fourier Transform (DFT) $F_n:\,\C_n(\Z) \to \C^{2n+1}$,
\[
(F_n [x])_k = \frac{1}{\sqrt{2n+1}} \sum_{ |\tau| \leq n} \exp\left(\frac{2\pi i k \tau }{2n+1}\right) x_\tau, \quad  1\leq k\leq 2n+1.
\]
While the new bounds are inflated (the additional factor $\rho$ is present in both bounds), the bound~\eqref{eq:intro-autoconv-norm} is essentially stronger than its counterpart~$\| F_m[\phi^o]\|_{1} \le \rho$ one could extract from~\eqref{eq:l2-oracle-norm}.

Based on this observation, the authors studied in~\cite{goldenshluger1997adaptive,jn1-2009,harchaoui2015adaptive} a class of adaptive convolution-type ``uniform-fit'' estimators
which correspond to filters obtained by minimizing the uniform norm of the Fourier-domain residual~$F_{n}[y - y * \vphi]$ constrained (or penalized) by the $\ell_1$-norm of the DFT of the filter.
%~$\|\F_{2m}[\vphi]\|_{1}$.
Such estimators can be efficiently computed since the corresponding filters are given as optimal solutions to well-structured convex optimization problems.

As it is common in adaptive nonparametric estimation, one can measure the quality of an adaptive estimator with the factor---the ``cost of adaptation''---by which the risk of such an estimator is greater than that of the corresponding ``oracle'' estimator which the adaptive one is trying to ``mimic''.
As it turns out, ``uniform-fit'' estimators studied in~\cite{goldenshluger1997adaptive,jn1-2009,harchaoui2015adaptive} admit the pointwise risk bounds similar to \rf{eq:ss1}, with extra factor~$C\rho^3\sqrt{\log (m+n)}$ as compared to \rf{eq:ss1} (see~\cite[Theorem~5]{harchaoui2015adaptive}).
On the other hand, there is a lower bound stating that the adaptation factor cannot be less than $c\rho\sqrt{\log m}$ when $m \ge c'n$~(cf.~\cite[Theorem~2]{harchaoui2015adaptive}), leaving the gap between these two bounds which may be quite significant when $\rho$ is large.
Furthermore, the choice of optimization objective (uniform fit of the Fourier-domain residual) in such estimators was dictated by the technical consideration allowing simpler control of the pointwise risk and seems artificial when the estimation performance is measured by the~$\ell_2$-loss. %
\paragraph{Contributions.}
In this paper, we propose a new family of adaptive convolution-type estimators.
These estimators utilize an adaptive filter which is obtained by minimizing the $\ell_2$-norm of the residual constrained or penalized by the $\ell_1$-norm of the DFT of the filter.
Similarly to uniform-fit estimators, new estimators can be efficiently computed via convex optimization routines.
%moreover, their formulation is easily amenable to first-order proximal algorithms.
We prove oracle inequalities for the $\ell_2$-loss of these estimators, which lead to the improved risk bounds compared to the case of uniform-fit estimators.
Note that signal simplicity, as per Definition~\ref{def:l2simple}, involves a special sort of time-invariance of the oracle estimate:
filter $\phi^o \in \C_m(\Z)$ in Definition \ref{def:l2simple} is assumed to be ``good'' (cf.~\rf{eq:l2-oracle-bias}) uniformly over  $|t|\leq m+n$, what  can be understood as some kind of ``approximate local shift-invariance'' of the signal to be recovered. In fact, this property of the signal is operational when deriving corresponding risk bounds for adaptive recoveries.
In the present paper, in order to derive the oracle inequalities we replace the assumption of signal simplicity, as per Definition~\ref{def:l2simple}, with an explicit {\em approximate (local) shift-invariance} (ASI) assumption.
In a nutshell, the new assumption states that the unknown signal admits (locally) a decomposition
$
x=x^{\cS}+\varepsilon
$
where~$x^{\cS}$ belongs to an  {\em unknown}  shift-invariant linear subspace $\cS \subset \C(\Z)$ of a small dimension, and the residual component~$\varepsilon$ is small in $\ell_2$-norm or~$\ell_\infty$-norm.
The remainder terms in the established oracle inequalities explicitly depend on the subspace dimension $s = \dim(\S)$ and the magnitude $\varkappa$ of residual component.
We also study the relationship between our ASI assumption and the notion of signal simplicity introduced in \cite{jn1-2009}: %(and formalized here in Definition~\ref{def:l2simple}):
\begin{itemize}
\item
On one hand,  approximately shift-invariant signals constitute a subclass of simple signals (in fact, the widest known to us such subclass to date).
In particular, a ``uniform'' version of ASI assumption, in which the residual component~$\varepsilon$ is bounded in $\ell_\infty$-norm, implies signal simplicity (cf. Definition \ref{def:l2simple}) with simple dependence of parameters $\rho$ and $ \theta$ of the class on the ASI parameters~$s$ and $\varkappa$.
This, in turn, allows to derive improved bounds for the pointwise and $\ell_2$-loss of novel adaptive estimators over the class of signals satisfying the ``uniform'' version of ASI assumption.
\item
On the other hand, all known to us examples of simple signals in $\C(\Z)$ are those of signals close to solutions of low-order linear homogeneous difference equations, see~\cite{jn2-2010}; such signals are close to small-dimensional shift-invariant subspaces.
New bounds on the $\ell_2$- and pointwise risk for such signals established in this work improve significantly over the analogous bounds for such signals obtained in \cite{jn2-2010,harchaoui2015adaptive}.
\end{itemize}

As an illustration, we consider an application of the proposed approach to the problem of denoising a {\em harmonic oscillation}---a sum of complex sinusoids with arbitrary (unknown) frequencies.
The known approaches~\cite{recht1,recht2} to this problem are based on the ideas from sparse recovery \cite{duarte2013spectral}
and impose {\em frequency separation} conditions to obtain sharp statistical guarantees (see Section~\ref{sec:sines} for more details).
In contrast, deriving near-optimal statistical guarantees for adaptive convolution-type estimators in this problem does not require this type of assumptions.

Preliminary versions of some results presented in this paper were announced in~\cite{ostrovsky2016structure}.

\paragraph{Manuscript organization.}
We present the problem of adaptive interpolation and prediction and  introduce necessary notation in Section~\ref{sec:main}.
In Section~\ref{sec:oracle2} we introduce adaptive estimators and present oracle inequalities for their $\ell_2$-loss. Then we use these inequalities to derive guarantees for $\ell_2$- and pointwise risks of adaptive estimates in Section \ref{sec:shift-inv}. In particular, in  Section \ref{sec:discuss-shift} we discuss the structure of the classes of approximately shift-invariant signals over $\Z$ and show that such signals are close, in certain sense, to complex exponential polynomials---solutions to linear homogeneous difference equations. %homogeneous linear difference equations.
We then specify statistical guarantees for adaptive interpolation and prediction of such signals; in particular, we establish new bounds for adaptive prediction of {\em generalized harmonic oscillations}
which are sums of complex sinusoids modulated by polynomials.
Finally, in Section \ref{sec:sines} we consider an application of the proposed estimates to the problem of full recovery of a generalized (or usual) harmonic oscillation, and compare our approach against the state of the art for this problem. %in~\cite{recht2}.
%and compare our approach against the state-of-the-art.
To streamline the presentation we defer technical proofs to appendix.

\section{Problem description}
\label{sec:main}
\subsection{Notation}
\label{sec:intro-notation}

%We now introduce the notation used in the subsequent exposition.
%Additional notation necessary for the proofs is presented in the appendix.

%\paragraph{Spaces.}
%$\Z$, $\R$, and $\C$ are the sets of all integer, real, and complex numbers correspondingly.
%$\Z^{+}$ is the set of all non-negative integers, $\Z^{++}$ the set of all strictly positive integers, idem for $\R$.
%$\R^n$ and $\C^n$ for $n \in \Z^+$ are the real and complex coordinate spaces.
%$\R^{m \times n}$ and $\C^{m \times n}$ are the spaces of $m \times n$ matrices with real or complex elements.

%\paragraph{Probability.}
%%$\Prob(E)$ and $\P(E)$ are interchangeably used for the probability of an event $E$.
%$\Prob(E)$ denotes the probability of an event $E$.
%When described by a predicate $\cC$, an event is written as $\{\cC\}$. The expectation is denoted by $\E$, and the probability measure over which it is taken is clear from context.
%%$\chi_n^2$ is the chi-square law with $n$ degrees of freedom.

%
%
%\paragraph{Matrices, vectors, and signal slices.}
We follow the ``Matlab convention'' for matrices: $[A, B]$ and $[A; B]$ denote, respectively, the horizontal and vertical concatenations of two matrices of compatible dimensions.
Unless explicitly stated otherwise, all vectors are column vectors.
Given  a signal $x \in \C(\Z)$ and $n_1, n_2 \in \Z$ such that $n_1 \le n_2$, we define the ``slicing'' map

\beq{intro-slicing}
x_{n_1}^{n_2} := [x_{n_1}; ...; x_{n_2}]. % \in \C^{n_2-n_1+1}.
\eeq
In what follows, when it is unambiguous, we use the shorthand notation $\tau\leq n$ ($\tau<n$, $|\tau|\leq n$, etc.) for the set of integers satisfying the inequality in question.
\paragraph{Convolution and filters.}
Recall that $\C(\Z)$ is the linear space of all two-sided complex sequences, and $\C_n(\Z)$ denotes the space of such sequences which vanish outside $[-n,...,n]$. We call the smallest $m \in \Z_+$ such that $\phi \in \C_m(\Z) $ the \textit{width} of $\phi$ and denote it $\Win(\phi)$.
Note that \eref{intro-slicing} allows to identify $\C_n(\Z)$, %$\C_m^+(\Z)$, $\C_m^-(Z)$, $\C_m^h(\Z)$
with complex vector space $\C^{2n+1}$. It is also convenient to identify $x\in \C(\Z)$ with its Laurent series
$
x(z)=\sum_{j} x_j z^j.
$
The (discrete) convolution of $\vphi * \psi \in \C(\Z)$ of $\vphi, \psi \in \C(\Z)$ is defined as
\[
[\vphi * \psi]_t := \sum_{\tau \in \mathbb{Z}} \vphi_\tau \psi_{t-\tau} %\quad t \in \Z.
\]
and is, clearly, a commutative operation.
One has $[\vphi * \psi](z)=\vphi(z) \psi(z)$ with
\[
\Win(\vphi * \psi) \le \Win(\vphi) + \Win(\psi).
\]
In what  follows, $\Delta$ stands for the forward shift operator on $\C(\Z)$:
\[
[\Delta x ]_t = x_{t-1},
\]
and $\Delta^{-1}$ for its inverse, the backward shift. Then
\[\vphi * \psi = \vphi(\Delta) \psi.
\]
%
%\paragraph{Convolution-type estimates.}
Given $\vphi\in \C(\Z)$ with $\Win(\vphi) < \infty$ and observations $y = (y_\tau)$, we can associate with $\vphi$ the linear estimate $\wh{x}$ of $x \in \C(\Z)$ of the form
\be
\wh{x} = \vphi * y = \vphi(\Delta) y
\ee{intest}
($\wh x$ is simply a kernel estimate over the grid $\Z$ corresponding to a finitely supported discrete kernel~$\vphi$).
% the linear estimate $[\vphi*y]_t$ with a finite-domain $\vphi$ may be alternatively expressed as $[\vphi(\Delta)y]_t$.
The just defined ``convolution'' (kernel) estimates are referred to as linear filters in signal processing; with some terminology abuse, we also call {\em filters} elements of $\C(\Z)$ with finitely many nonzero entries.

\paragraph{Norms.} For $x,\,y\in \C(Z)$ we denote $\langle x, y \rangle$ the Hermitian inner product
$ \langle x, y\rangle= \sum_{\tau \in \Z} \overline{x}_\tau y_\tau$,
$\overline{x}_\tau$ being the complex conjugate of $x_\tau$; for $n\in \Z_+$ we put
\[
\langle x, y \rangle_n = \sum_{|\tau|\leq n} \overline{x}_\tau y_\tau.
\]
Given $p \ge 1$ and~$n \in \Z_+$ we define  semi-norms on $\C(\Z)$ as follows:
\[
\|x\|_{n,p} := \left(\sum_{|\tau|\leq n} | x_\tau |^p\right)^{1/p}
\]
with $\|x\|_{n,\infty}=\max_{|\tau|\leq n} |x_\tau|$.
When such notation is unambiguous, we also use $\|\cdot\|_p$ to denote the ``usual'' $\ell_p$-norm on $\C(\Z)$, e.g.,~$\|x\|_p=\|x\|_{n,p}$ whenever $\Win(x) \le n$.
\par
We define the (unitary) Discrete Fourier Transform (DFT) operator $F_n: \C_n(\Z) \to \C^{2n+1}$ by %{\cbl Verify if $F_n$ is $\C_n(\Z) \to \C_n(\Z)$ of $\C_n(\Z) \to \C^{2n+1}$ }
\[
\left(F_n [x]\right)_k = \frac{1}{\sqrt{2n+1}} \sum_{|\tau| \leq n} \exp\left(-\frac{i2\pi k \tau }{2n+1}\right) x_\tau, \quad 1\leq k \leq 2n+1.
\]
The unitarity of DFT implies
%that $F_n[\zeta]$ follows the same law as~$[\zeta]_{-n}^n$, and also
the Parseval identities: for any $x, y \in \C(\Z)$ and $n \in \Z_+$ one has
\begin{equation}
\label{eq:intro-parseval}
\langle x, x \rangle_n = \langle F_n [x], F_n [x] \rangle, \quad
\|x\|_{n,2} = \|F_n [x]\|_{2}.
\end{equation}
In what follows, $c$, $C$, $C'$, etc., stand for absolute constants whose exact values can be recovered from the proofs. We use the~$O(\cdot)$ notation: for two functions $f,g$ of the same argument~$t$, $f = \cO(g)$ means that there exists $C<\infty$ such that $|f(t)| \le C|g(t)|$ for all~$t$ in the domain of $f$.
%Besides this, $f = \cO_\ell(g)$ means that $f \le g$ holds up to a logarithmic factor in the common argument of $f$ and $g$.

\subsection{Problem statement}
\label{sec:problem-statement}
We consider the problem of estimating the signal $x\in \C(\Z)$ given noisy observations $y_\tau := x_\tau + \sigma \zeta_\tau$ on the segment $-L\leq \tau\leq L$ (cf. \rf{eq:intro-observations});
here $\zeta_t \sim \CN(0,1)$ are i.i.d.~standard complex-valued Gaussian random variables.  Here we discuss different settings of this problem:
\begin{itemize}
\item {\em Signal interpolation} in which, when computing the estimate of $x_t$, one can use observations {both on the left and on the right} of $t$.
For the sake of simplicity, we consider the ``symmetric'' version of this problem where the objective is, given $|m|\leq L$, to build an estimate $\wh x_t=[\wh \vphi * y]_t$ of $x_t$ for~$|t|\leq L-m$, with $\wh \vphi\in \C_m(\Z)$ depending on observations.
\item {\em Signal prediction} in which, when computing the estimate of $x_t$, we are allowed to use observations only on one side of $t$, e.g.,~observations for $\tau\leq t-h$ where $h\in \Z_+$ is a given {\em prediction horizon}. For the sake of clarity, in this paper we only consider the version of this problem with $h=0$ (often referred as {\em filtering} in signal processing literature); the general situation can be treated in the same way at the expense of more involved notation. In other words, we are looking to build a data-driven filter $\wh \vphi\in \C_m(\Z)$ and the ``left'' estimate of $x_t$, $-L+2m\leq t\leq L$ (utilizing observations $y_\tau,\,\tau\leq t$),
    \[
    \wh x_t=\sum_{\tau=0}^{2m}\vphi_{\tau-m} y_{t-\tau}=\sum_{s=-m}^{m}\vphi_{s} y_{t-s-m}=[\vphi * (\Delta^{m}y)]_t.
    \]
The corresponding ``right'' estimate of $x_t$, $-L\leq t\leq L-2m$ (utilizing observations $y_\tau,\,\tau\geq t$) writes
\[
    \wh x_t=\sum_{\tau=0}^{2m}\vphi_{m-\tau} y_{t+\tau}=\sum_{s=-m}^{m}\vphi_{s} y_{t-s+m}=[\vphi * (\Delta^{-m}y)]_t.
    \]
\end{itemize}
Given a set $\X$ of signals, $m, n\in\Z_+$, observations $y_\tau$ for $|\tau|\leq L=m+n$,  and the target estimation domain $D_n$ of length $2n+1$ (e.g.,~$D_n=\{-n,...,n\}$ in the case of signal interpolation, or $D_n=\{-n+m,...,n+m\}$ in the case of filtering), we quantify the accuracy of estimate $\wh x$ using two types of risks:
\begin{itemize}
\item{\em maximal over $\X$  $\ell_2$ (integral) $\alpha$-risk:} the smallest maximal over $x\in \X$ radius of $(1-\alpha)$-confidence ball of $\|\cdot\|_2$-norm on $D_n$ centered at $\wh x$:
 \[%bse
\risk_{D_n,2,\alpha}(\wh x|\X)=\inf \left\{r:\sup_{x\in \X}\Prob\left\{\left(\sum_{t\in D_n} |[\wh x-x]_t|^2\right)^{1/2}\geq r\right\}\leq \alpha\right\};
\]
%ese
\item{\em maximal over $\X$ pointwise $\alpha$-risk:}  the smallest maximal over $x\in \X$ and $t\in D_n$ $(1-\alpha)$-confidence interval for $x_t$ centered at $\wh x_t$:
\[%bse
\risk_{D_n,\alpha}(\wh x|\X)=\inf \left\{r:\sup_{x\in \X}\Prob\left\{|[\wh x-x]_t|\geq r\right\}\leq \alpha\, \; \forall\, t\in D_n\right\}.
%\\\risk_{\epsilon}(\wh x|X)&=&\inf \left\{r:\sup_{x_*\in X}\Prob\{\|\wh x-x_*\|\geq r\}\leq \epsilon\right\}.
\]
When $n=0$ the estimation interval  $D_n=\{t\}$ is a singleton, and  the latter definition becomes that of the ``usual'' worst-case over $\X$ $(1-\alpha)$-confidence interval for $x_t$:
\[%bse
\risk_{\alpha}(\wh x_t|\X)=\inf \left\{r:\sup_{x\in \X}\Prob\left\{|[\wh x-x]_t|\geq r\right\}\leq \alpha\right\}.
%\\\risk_{\epsilon}(\wh x|X)&=&\inf \left\{r:\sup_{x_*\in X}\Prob\{\|\wh x-x_*\|\geq r\}\leq \epsilon\right\}.
\]

\end{itemize}

%\paragraph{Preliminaries.}
\section{Oracle inequalities for $\ell_2$-loss of adaptive estimators}
\label{sec:oracle2}
\subsection{Adaptive signal interpolation}
\subsubsection{Adaptive recoveries}
Given $m,n \in \Z_+$, $L=m+n$, and $\overline\vrho > 0$, %let $\wh\vphi_{\con}\C_m(\Z)$
consider the optimization problem
%an optimal solution to the following optimization problem:
\begin{equation}
\label{opt:l2con}
\min_{\vphi\in \C_{m}(\Z)}%\left\{
\|y-\varphi*y\|_{n,2}^2 \;\;
\mbox{subject to}\;\;
%\phi\in \Phi(\bar\vrho)=\left\{
%:\;
\|F_m[\varphi]\|_{1} \leq \frac{\bar{\varrho}}{\sqrt{2m+1}}.%\right\}
\tag{{\bf Con}}
\end{equation}
Note that \rf{opt:l2con} is clearly solvable; we denote $\wh\vphi_{\con}$ its optimal solution and refer to
\[\widehat x_{\con} = \wh\vphi_{\con} * y\]
as the {\em constrained (least-squares) estimate} of $x$.
%In the sequel, we will prove a sharp oracle inequality for the constrained estimator, which states that the $\ell_2$-loss of this estimator is comparable to the $\ell_2$-loss of \textit{any} filter $\varphi$ feasible to~\eqref{opt:l2con}, and in particular, to~$\vphi^o$ provided that $\bar \vrho = \vrho$.
Computing $\wh\vphi_{\con}$ requires setting the problem parameter $\bar\vrho$ which, ideally, would be set proportional to the $\ell_1$-norm of the DFT of some ideal (oracle) filter, or a non-trivial upper bound on it.
Because this is not often possible in practice, we also consider the {\em penalized estimator}
\[
\wh x_\pen = \wh\vphi_{\pen} * y,
\]
where, for $\lambda >0$, $\wh\vphi_{\pen}\in \C_m(\Z)$ is selected as an optimal solution to the (solvable) problem
\begin{equation}
\label{opt:l2pen++}
\min_{\varphi \in \C_{m}(\Z) }
\|y-\varphi*y\|^{2}_{n,2} + \sigma^2 \lambda^2  (2m+1) \|F_m[\varphi]\|_{1}^2  \tag{{\bf Pen}}.
\end{equation}
Instead of knowing $\bar\vrho$, some knowledge of noise variance~$\sigma^2$ is required to tune this estimator.
%Below, we show that as long the penalization parameter $\lambda$ is set as an absolute constant, the $\ell_2$-loss of the penalized estimator enjoys essentially the same bound as the constrained estimator with the choice of $\bar{\vrho}$ corresponding to the optimal bias-variance balance of the oracle estimate.
Hence, the practical recommendation is to use~\eqref{opt:l2pen++} when %possible, \ie~whenever~
$\sigma^2$ is known or can be estimated.

\subsubsection{Oracle inequalities for $\ell_2$-loss }\label{sec:l2oracle}
Despite striking similarity with Lasso estimators \cite{tibshirani,candes2007danzig,bickel2009simultaneous}, the proposed estimates are of quite different nature.
First of all, solving optimization problems \rf{opt:l2con} and \rf{opt:l2pen++} allows to recover a filter but not the signal itself, and this filter is generally not sparse neither in time nor in Fourier domain (unless the signal to recover is a sum of harmonic oscillations with frequencies on the ``DFT grid'').
Second, the equivalent of ``regression matrices'' involved in these procedures cannot be assumed to satisfy any kind of ``restricted incoherence'' conditions usually imposed to prove statistical properties of ``classical'' $\ell_1$-recovery routines (see~\cite[Chapter 6]{buhlmann2011statistics} for a comprehensive overview of such conditions).
Moreover, being constructed from {\em noisy observations}, these matrices depend on the noise, which poses some extra difficulties in the analysis of the adaptive estimates, in particular, leading to the necessity of imposing some restrictions on the signal class.

In what follows, when analyzing adaptive estimators we constrain the unknown signal $x$ on the interval $|\tau|\leq L$ to be ``close'' to some shift-invariant linear subspace $\S$. Specifically, consider the following assumption:

\begin{ass}[{Approximate local shift-invariance}]
\label{ass:subspace}
We suppose that $x \in \C(\Z)$ admits a decomposition
\[
x=x^{\S}+\varepsilon.
\]
Here, $ x^{\S}\in \S$ where $\S$ is some (unknown) shift-invariant linear subspace of $\C(\Z)$ with $s:= \dim(\S) \le 2n+1$, and $\varepsilon$ is bounded in the $\ell_2$-norm: for some $\varkappa \ge 0$ one has
\begin{align}\label{eq:eps-small}
\left\|\Delta^{-\tau} \varepsilon\right\|_{n,2} \leq \varkappa \sigma, \quad |\tau|\leq m.
\end{align}
We denote $\X_{m,n}(s,\varkappa)$ the class of such signals.
\end{ass}
\paragraph{Remarks.} Assumption~\ref{ass:subspace} merits some comments.
\par
Observe that $\X_{m,n}(s,\varkappa)$ is in fact the subset of $\C(\Z)$ comprising sequences which are close, in the sense of~\rf{eq:eps-small}, to {\em all} $s$-dimensional shift-invariant subspaces of $\C(\Z)$. Similarly to Assumption~\ref{ass:subspace}, signal ``simplicity'' as set by Definition \ref{def:l2simple} also postulates a kind of ``local time-invariance'' of the signal: it states that there exists a linear time-invariant filter which reproduces the signal ``well'' on a certain interval. However, the actual relationship between the two notions is rather intricate and will be discussed in Section~\ref{sec:shift-inv}.

Letting the signal to be close, in $\ell_2$-norm, to a shift-invariant subspace---instead of simply belonging to the subspace---extends the set of signals and allows to address nonparametric situations. As an example, consider discretizations over a uniform grid in $[0,1]$ of functions from the Sobolev ball. Locally, such signals are close in $\ell_2$-norm to polynomials on the grid which satisfy a linear homogeneous difference equation and hence belong to a shift-invariant subspace of small dimension \cite{jn2-2010}.
%More generally, Assumption \ref{ass:subspace} holds with~$\varkappa = 0$ for
Other classes of signals for which  Assumption~\ref{ass:subspace} holds are discretizations of complex sinusoids modulated with smooth functions and signals satisfying linear difference inequalities \cite{jn2-2010}.

We now present \textit{oracle inequalities} which relate the $\ell_2$-loss of adaptive filter $\wh\varphi$ with the best loss of any feasible solution~$\varphi$ to the corresponding optimization problem.
These inequalities, interesting for their own sake, are also operational when deriving bounds for the  pointwise and $\ell_2$-losses of the proposed estimators.
We first state the result for the constrained estimator.
\begin{theorem}
\label{th:l2con}
Let $s, m,n \in \Z_+$, $\varkappa\ge 0$. Suppose that $x\in \X_{m,n}(s, \varkappa)$ and $\vphi$ is feasible for ~\eqref{opt:l2con}. Let $\wh\vphi_{\con}$ be an optimal solution to~\eqref{opt:l2con} with some $\bar\vrho > 1$, and let $\wh x_\con=\wh\vphi_{\con} *y$.
%so that
%\[
%\|\vphi\|^{\F}_{m,1}\leq \frac{\bar{\varrho}}{\sqrt{2m+1}}.
%\]
Then for any $\alpha\in ]0,1[$
% $0 < \alpha \le 1$  there exists set $\Xi_\alpha$ of noise realizations $\zeta$ such that $\Prob\{\Xi_\alpha\} \ge 1-\alpha$, and whenever $\zeta \in \Xi_\alpha$ the estimate
 it holds with probability at least $1-\alpha$: %$\wh x_{\con}$ satisfies
\be
\left\|x -  \wh x_{\con}\right\|_{n,2} &\le& \left\|x-\vphi * y\right\|_{n,2}\nn&&
	+ C\sigma\Big(\bar\vrho (\kappa^2_{m,n}+1)\log[(m+n)/\alpha] + \bar\vrho \varkappa\sqrt{\log[1/\alpha]} + s\Big)^{1/2}
\ee{eq:th1}
where
\[
\kappa_{m,n} := \sqrt\frac{2n+1}{2m+1}.
\]
\end{theorem}

The counterpart of Theorem~\ref{th:l2con} for the penalized estimator is as follows.

\begin{theorem}\label{th:l2pen++}
Let $s, m,n \in \Z_+$, $\varkappa,\lambda>0$. Suppose that $x\in \X_{m,n}(s, \varkappa)$ and $\vphi\in \C_m(\Z)$ with $\vrho(\vphi) = \sqrt{2m+1} \|F_m[\vphi]\|_{1}$. Let $\wh\vphi_{\pen}$ be an optimal solution to~\eqref{opt:l2pen++}.
Then for any $\alpha\in \,]0,1[$
% $0 < \alpha \le 1$  there exists set $\Xi_\alpha$ of noise realizations $\zeta$ such that $\Prob\{\Xi_\alpha\} \ge 1-\alpha$, and whenever $\zeta \in \Xi_\alpha$ the estimate
the estimate $\wh x_\pen=\wh\vphi_{\pen} *y$ satisfies  with probability at least $1-\alpha$:
\begin{equation}\label{eq:l2pen++}
\left\|x - \wh x_{\pen}\right\|_{n,2} \le
\|x-\vphi * y\|_{n,2} + \sigma \Big( \lambda \vrho(\vphi) + C_1  \Rem_1/\lambda + C_2\Rem_2^{1/2}(\vphi)\Big)
\end{equation}
where
\begin{equation}
\label{eq:remainder-pen}
\begin{aligned}
\Rem_1 = \Rem_1(\varkappa,\kappa_{m,n},\alpha) &= (\kappa_{m,n}^2+1)\log[(m+n)/\alpha] + \varkappa\sqrt{\log[1/\alpha]} + 1, \\
\Rem_2(\vphi) = \Rem_2(\vphi,s,\varkappa,\alpha) &= \varrho(\vphi)\log[1/{\alpha}] + \varkappa\sqrt{\log[1/\alpha]} + s.
\end{aligned}
\end{equation}
In particular, when setting $\lambda = \Rem_1^{1/2}$ we obtain
\begin{align*}
\left\|x - \wh x_\pen\right\|_{n,2}
&\le \|x-\vphi * y\|_{n,2} + C\sigma\Big(\Rem_1^{1/2}\varrho(\vphi) + \Rem_2^{1/2}(\vphi)\Big).
%&\le \|x-\vphi * y\|_{n,2} + C'\sigma \Rem_0[\vrho^2]^{1/2}.
\end{align*}
%Note that with $\Rem_0[r]$ defined in~\eqref{eq:remainder-con}.
\end{theorem}
One may observe that, ideally, $\bar\vrho$ in \rf{opt:l2con} should be selected as $\vrho(\vphi^o)=\sqrt{2m+1}\|F_m[\vphi^o]\|_1$ where $\vphi^o$ is an ideal ``oracle filter,'' while the penalty parameter in \rf{opt:l2pen++} would be set to $\lambda=[{C_1\Rem_1}/\vrho(\vphi^o)]^{1/2}$. These choices would result in the same remainder terms in
\rf{eq:th1} and \rf{eq:l2pen++}
(order of $\sigma(\vrho(\vphi^o)(1+\varkappa)+s)^{1/2}$ up to logarithmic factors). Obviously, this choice cannot be implemented since the value $\vrho(\vphi^o)$ is unknown.
Nevertheless, Theorem~\ref{th:l2pen++} provides us with an implementable choice of $\lambda$ that still results in an oracle inequality, at the expense of a larger remainder term which now scales as $\sigma[\vrho(\vphi^o)\sqrt{1+\varkappa} + \sqrt{s}]$.

%\section{Extension}\label{sec:extenstions}
\subsection{Adaptive signal filtering}
\label{sec:prediction}
Here we consider the ``left'' version of the problem in which we are given observations~$(y_\tau)$ on the interval $-L\leq \tau \le L$, and our objective is to build a (left) convolution estimate $\wh x_t=[\wh \vphi * (\Delta^{m}y)]_t$ of $x_t$, $t\in\{-L+2m\leq t\leq L\}$, using an observation-driven filter $\wh\vphi \in \C_{m}(\Z)$. Clearly, the treatment of the ``right'' version of the problem is completely analogous up to obvious modifications.
Let us consider the following counterparts of~\eqref{opt:l2con} and~\eqref{opt:l2pen++}: %modified for the prediction setting:
\begin{align}
\label{opt:l2con-pred}
&\min_{\vphi\in \C_m(\Z)}\left\| \Delta^{-m} (y-\varphi*\Delta^m y)\right\|_{n,2}^2 \mbox{ subject to }
\big\|F_m[\varphi]\big\|_{1} \leq \frac{\bar{\vrho}}{\sqrt{2m+1}}
\tag{\bf{Con}$^+$},\\
\label{opt:l2pen-pred++}
&\min_{\vphi\in \C_m(\Z)} \left\| \Delta^{-m} (y-\varphi*\Delta^m y)\right\|_{n,2}^2 +
\sigma^2 \lambda^2 (m+1) \big\|F_m[\varphi] \big\|_{1}^2
\tag{\bf{Pen}$^+$}.
%\end{equation}
\end{align}
Same as in the interpolation setting, both problems are clearly solvable, so their respective optimal solutions $\widehat\varphi_{\con}$ and $\widehat\varphi_{\pen}$ are well-defined.
A close inspection of the proofs of Theorems~\ref{th:l2con} and \ref{th:l2pen++} shows that their results remain valid, with obvious adjustments, in the setting of this section. Namely, we have the following analog of those statements.
\begin{prop}
\label{prop:l2all}
Let $s, m,n \in \Z_+$, $\varkappa\ge 0$, and $x\in \X_{m,n}(s, \varkappa)$; let $\alpha\in \,]0,1[$.
\begin{enumerate}
\item Let $\bar\vrho > 1$ be fixed, $\varphi$ be feasible to \rf{opt:l2con-pred}, and let  $\wh x_\con=\wh\vphi_{\con}*\Delta^{m}y$ where $\wh\vphi_{\con}$ is an optimal solution to~\eqref{opt:l2con-pred}; then with probability at least $1-\alpha$
 estimate
$\wh x_{\con}$ satisfies
\bse
\left\|\Delta^{-m}(x -  \wh x_{\con})\right\|_{n,2} &\le& \left\|\Delta^{-m}(x-\vphi * \Delta^my)\right\|_{n,2}\nn&& + C\sigma\Big(\bar\vrho (\kappa^2_{m,n}+1)\log[(m+n)/\alpha] + \bar\vrho \varkappa\sqrt{\log[1/\alpha]} + s\Big)^{1/2}.
\ese
\item Let $\vphi\in\C_m(\Z)$ with $\vrho(\vphi) = \sqrt{2m+1} \|F_m[\vphi]\|_{1}$, and let $\wh x_\pen=\wh\vphi_{\pen}*\Delta^{m}y$ where $\wh\vphi_{\pen}$ is an optimal solution to~\eqref{opt:l2pen-pred++} with $\lambda>0$; then
$\wh x_{\pen}$ satisfies with probability at least $1-\alpha$
\bse
\left\|\Delta^{-m}(x - \wh x_{\pen})\right\|_{n,2} &\le &
\|\Delta^{-m}(x-\vphi * \Delta^{m}y)\|_{n,2}\\&&+\sigma \Big( \lambda \vrho(\vphi) + C_1  \Rem_1/\lambda + C_2\Rem_2^{1/2}(\vphi)\Big)
\ese
where $\Rem_1$ and $\Rem_2(\vphi)$ are defined in \rf{eq:remainder-pen}.

\end{enumerate}
\end{prop}

\section{Risk bounds for adaptive recovery under ASI}
\label{sec:shift-inv}
In order to transform the oracle inequalities of Theorems \ref{th:l2con}, \ref{th:l2pen++} and Proposition \ref{prop:l2all} into risk bounds for adaptive recoveries, we need to establish bounds for oracle risks on the classes of approximately shift-invariant signals. We start with the interpolation setting.
\subsection{Risk bounds for adaptive signal interpolation}
Results of this section are direct corollaries of the following statement which may be of independent interest.
\begin{prop}
\label{prop:shift-inv-bilateral}
Let $\S$ be a shift-invariant subspace of $\C(\Z)$ of dimension~$s \le m+1$.
Then there exists a filter $\phi^o \in \C_m(\Z)$ such that for all~$x\in \S$ one has~$x=\phi^o * x$ and
\[
\|\phi^o\|_2\leq \sqrt{\frac{2s}{2m+1}}.
\]%ee{eq:interps}
\end{prop}
\noindent In other words, signals $x\in \S$ are $(m,n,\rho,0)$-simple in the sense of Definition \ref{def:l2simple}, for any~$n \in \Z_+$ and $m\geq s-1$, with $\rho = \sqrt{2s}$ and $\theta=0$.
\par
When combined with Theorems \ref{th:l2con} and \ref{th:l2pen++}, Proposition \ref{prop:shift-inv-bilateral} implies the following bound on the integral risk of adaptive recovery.
\begin{prop}
\label{pr:l2sloss}
Let $s, m,n \in \Z_+$, $m\geq 2s-1$, $\varkappa\ge 0$, and let $D_n=\{-n,...,n\}$.
\item[(i)] Assume that $\wh x_{\con}=\wh\vphi_{\con}*y$ where $\wh\vphi_{\con}$ is an optimal solution  to~\eqref{opt:l2con} with some $\bar\vrho \geq 4s$.
Then for any $\alpha\in ]0,1/2]$
\[
\risk_{{D_n},2,\alpha}(\wh x_{\con}|\X_{m,n}(s,\varkappa))\leq C\psi^{\alpha}_{m,n}(\sigma,s,\varkappa;\bar\vrho)
\]%ee{eq:riskl2-1}
where
\bse
\psi^{\alpha}_{m,n}(\sigma,s,\varkappa;\bar\vrho)&=&\sigma s\big(\kappa_{m,n}\sqrt{\log[1/\alpha]}+\varkappa\big)\\&&+\sigma\Big(\bar\vrho (\kappa^2_{m,n}+1)\log[(m+n)/\alpha] + \bar\vrho \varkappa\sqrt{\log[1/\alpha]} + s\Big)^{1/2}.
\ese
In particular, when $\bar\vrho\leq C's$ is chosen in~\eqref{opt:l2con} one obtains
\be
\risk_{{D_n},2,\alpha}(\wh x_{\con}|\X_{m,n}(s,\varkappa))\leq C\overline\psi^{\alpha}_{m,n}(\sigma,s,\varkappa)
\ee{eq:riskl2-2}
with
\bse
\overline \psi^{\alpha}_{m,n}(\sigma,s,\varkappa)&=&\sigma s\big(\kappa_{m,n}\sqrt{\log[1/\alpha]}+\varkappa\big)\\&&+\sigma\Big(s (\kappa^2_{m,n}+1)\log[(m+n)/\alpha] + s \varkappa\sqrt{\log[1/\alpha]} + s\Big)^{1/2}.
\ese
\item[(ii)]
Let $\lambda = \Rem_1^{1/2}$ with $\Rem_1$ as defined in \rf{eq:remainder-pen}, and let $\wh x_{\pen}=\wh\vphi_{\pen}*y$ where $\wh\vphi_{\pen}$ is an optimal solution  to~\eqref{opt:l2pen++}.
Then for any $\alpha\in (0,1/2]$
\[
\risk_{{D_n},2,\alpha}(\wh x_{\pen}|\X_{m,n}(s,\varkappa))\leq C\widetilde\psi^{\alpha}_{m,n}(\sigma,s,\varkappa)
\]%{eq:riskl2-3}
where
\bse
\widetilde\psi^{\alpha}_{m,n}(\sigma,s,\varkappa)&=&\sigma s\big(\kappa_{m,n}\sqrt{\log[1/\alpha]}+\varkappa\big)
+\sigma s(\kappa_{m,n}+1)\sqrt{\log[(m+n)/\alpha]}.\ese
\end{prop}
We are now ready to derive  bounds for the pointwise risk of adaptive estimates described in the previous section. To establish such bounds we need to replace Assumption \ref{ass:subspace} with a somewhat stronger uniform analog.
\begin{ass}[{Approximate locally {\em uniform} shift-invariance}]
\label{ass:subspace_infty}
Let $n\geq m\in \Z_+$. We suppose that $x \in \C(\Z)$ admits a decomposition
\[
x=x^{\S}+\varepsilon.
\]
Here $ x^{\S}\in \S$ where $\S$ is some (unknown) shift-invariant linear subspace of $\C(\Z)$ with $s:= \dim(\S) \le 2n+1$, and $\varepsilon$ is {\em uniformly} bounded: for some $\varkappa \ge 0$ one has
\begin{align}
\label{eq:eps-small-infty}
|\varepsilon_\tau| \leq {\varkappa \sigma\over \sqrt{2n+1}}, \quad |\tau|\leq n+m.
\end{align}
We denote $\overline \X_{m,n}(s,\varkappa)$ the class of such signals.
\end{ass}
\noindent Observe that if $x\in \overline \X_{m,n}(s,\varkappa)$ then also $x\in \X_{m,n}(s,\varkappa)$. Therefore, the bounds of Proposition \ref{prop:shift-inv-bilateral} also hold true for the risk of adaptive recovery on $\overline \X_{m,n}(s,\varkappa)$. Furthermore, bound \rf{eq:eps-small-infty} of Assumption \ref{ass:subspace_infty} now leads to the following bounds for pointwise risk of recoveries $\wh x_\con$ and $\wh x_{\pen}$.
\begin{prop}
\label{pr:linfsloss}
Let $s, m,n \in \Z_+$ with  $m\geq 2s-1$ and $n\geq \lfloor m/2\rfloor$ (here $\lfloor \cdot\rfloor$ stands for the integer part),   $\varkappa\ge 0$; let also $D_{n,m}=\{-n+\lfloor m/2\rfloor,...,n-\lfloor m/2\rfloor\}$.
\item[(i)] Let $\wh x_{\con}=\wh\vphi_{\con}*y$ where $\wh\vphi_{\con}$ is an optimal solution  to~\eqref{opt:l2con} with $\bar\vrho \in[4s, Cs]$ for some $C\geq 4$.\footnote{For the sake of conciseness, here we only present the result for the constrained recovery with $\bar \vrho \asymp s$.}
Then for any $\alpha\in ]0,1/2]$
\be
\risk_{{D_{n,m}},\alpha}(\wh x_{\con}|\overline\X_{m,n}(s,\varkappa))\leq C'\overline\varsigma^\alpha_{m,n}(\sigma,s,\varkappa)
\ee{eq:pi}
where
\bse \overline\varsigma^\alpha_{m,n}(\sigma,s,\varkappa)\hspn&=&\hspn{\sqrt{s\over 2m+1}}\overline\psi^{\alpha}_{m,n}(\sigma,s,\varkappa)+{s\sigma\over \sqrt{2m+1}}\left(
\sqrt{s}\varkappa+\sqrt{\log \left[{(2m+1)/\alpha}\right]}+\sqrt{s\log[1/\alpha]%\left[{1\over \alpha}\right]
}\right)\\
\hspn&\leq&\hspn C''{s\sigma\over \sqrt{2m+1}}\left(\kappa_{m,n}\sqrt{s\log\left[{1/\alpha}\right]}+\varkappa+\kappa_{m,n}\sqrt{\log\left[{(m+n)/ \alpha}\right]}\right).
\ese
\item[(ii)]
Let $\wh x_{\pen}=\wh\vphi_{\pen}*y$ where $\wh\vphi_{\pen}$ is an optimal solution  to~\eqref{opt:l2pen++} with $\lambda = \Rem_1^{1/2}$, $\Rem_1^{\vphantom{1/2}}$ being defined in \rf{eq:remainder-pen}.
Then for any $\alpha\in (0,1/2]$
\[\risk_{{D_{n,m}},\alpha}(\wh x_{\pen}|{\overline\X_{m,n}(s,\varkappa)})\leq C\widetilde\varsigma^\alpha_{m,n}(\sigma,s,\varkappa)
\]
where
\bse
\widetilde\varsigma^\alpha_{m,n}(\sigma,s,\varkappa) \hspn&=&\hspn {\sqrt{s\over 2m+1}}\widetilde\psi^{\alpha}_{m,n}(\sigma,s,\varkappa)+{s\sigma\over \sqrt{2m+1}}\left(
\sqrt{s}\varkappa+\sqrt{\log \left[{(2m+1)/\alpha}\right]}+\sqrt{s\log\left[{1/\alpha}\right]}\right)\\
\hspn &\leq& \hspn {C's\sigma\over \sqrt{2m+1}}\left(\sqrt{s}\left(\kappa_{m,n}\sqrt{\log\left[{1/\alpha}\right]}
+\varkappa+\sqrt{\varkappa\log\left[{1/\alpha}\right]}\right)+\kappa_{m,n}\sqrt{\log\left[{(m+n)/\alpha}\right]}\right).
\ese
\end{prop}
\paragraph{Remark.} The above bounds for the pointwise risk of adaptive estimates may be compared against available lower bound and bounds for the risk of the uniform-fit adaptive estimate in the case where the signal to recover is a sum of $s$ {complex sinusoids}. In this situation, \cite[Theorem 2]{harchaoui2015adaptive} states the lower bound
$
c\sigma s \sqrt{\frac{\log m}{m}}$ for the pointwise risk of estimation with the upper bound
\[
\cO\left(\sigma s^3 \log[s] \sqrt{\frac{\log m}{m}}\right)
\] up to a logarithmic in $\alpha$ factor (cf. \cite[Section 4]{harchaoui2015adaptive}). Because the signal in question belongs to a $2s$-dimensional shift-invariant subspace of $\C(\Z)$, the {bound on} the pointwise risk in Proposition \ref{pr:linfsloss} results (recall that we are in the situation of $\varkappa=0$) in the bound
\[
\cO\left({\sigma s {\sqrt{s+\log m \over m}}}\right)\]
for adaptive estimates $\wh x_{\con}$ and $\wh x_{\pen}$ with significantly improved dependence on $s$.
\subsection{Risk bounds for adaptive signal filtering}
\label{sec:discuss-shift}
Our next {goal} is to bound the risk of the  constrained and penalized adaptive filters. Recall that in order to obtain the corresponding bounds in the interpolation setting we first established the result of Proposition \ref{prop:shift-inv-bilateral} which allows to bound the error of the oracle filter on any $s$-dimensional shift-invariant subspace of $\C(\Z)$. This result, along with oracle inequalities of Theorems \ref{th:l2con} and \ref{th:l2pen++}, directly led us to the bounds for the risk of adaptive interpolation estimates.
In order to reproduce the derivation in the previous section we first need to establish a fact similar to  Proposition \ref{prop:shift-inv-bilateral} which would guarantee existence of a {\em predictive} filter of small $\ell_2$-norm exactly reproducing all signals from any shift-invariant subspace of $\C(\Z)$. However, as we will see in an instant, the prediction case is rather different from the interpolation case: generally, a ``good  predictive filter'' one may look for---a reproducing predictive filter of small norm---simply does not exist in the case of prediction. And analysis of situations where such filter does exist is quite different from the simple  proof of  Proposition \ref{prop:shift-inv-bilateral}.
This is why, before returning to our original problem, it is useful to get a better understanding of the structure of shift-invariant subspaces of $\C(\Z)$.

\subsubsection{Characterizing  shift-invariant subspaces of $\C(\Z)$}
We start with the following
\begin{prop}
\label{th:shift_invariant} Solution set of a homogeneous linear difference equation
\begin{equation}\label{eq:shift-inv-diff_eq}
[p(\Delta) x]_t \left[= \sum_{\tau = 0}^s p_\tau x_{t-\tau} \right] = 0, \quad t \in \Z,
\end{equation}
with a characteristic polynomial $p(z) = 1 + p_1 z + ... + p_s z^s$ is a shift-invariant subspace of $\C(\Z)$ of dimension at most~$s$.
%\item
\\
Conversely, any shift-invariant subspace of $\C(\Z)$ of dimension $s$ is the solution set of a difference equation of the form~\eqref{eq:shift-inv-diff_eq} with $\deg(p)=s$; such polynomial is unique if normalized by~$p(0)=1$.
%\end{enumerate}
\end{prop}
\noindent Recall that the set of solutions of equation \eqref{eq:shift-inv-diff_eq} is spanned by \textit{exponential polynomials}. Namely,
let $z_k$, for $k=1,\, ...,\, r \le s$,  be the distinct roots of $p(z)$ with corresponding multiplicities $m_k$, and let~$\omega_k \in \C$ be such that $z_k = e^{-i \omega_k}$.
Then solutions to~\eqref{eq:shift-inv-diff_eq} are exactly sequences of the form
\[%begin{equation}
%\label{eq:exp-poly}
x_t = \sum_{k=1}^r q_k(t) e^{i \omega_k t}
\]%end{equation}
where $q_k(\cdot)$ are arbitrary polynomials of $\deg(q_k) = m_k-1$.
For instance, discrete-time polynomials of degree $s-1$ satisfy~\eqref{eq:shift-inv-diff_eq} with~$\textup{p}(z) = (1-z)^s$;
another example is that of~\textit{harmonic oscillations} with given (all distinct) $\omega_1,...,\omega_s\in [0,2\pi[$,
\begin{equation}
\label{eq:shift-inv-harmonic}
x_t = \sum_{k=1}^s q_k \e^{i \omega_k t}, \quad\quad q\in \C^s,
\end{equation}
which satisfy~\eqref{eq:shift-inv-diff_eq} with~$p(z) =\prod_{k=1}^{s} (1-e^{i \omega_k}z)$.
Thus, the set of complex harmonic oscillations with fixed frequencies $\omega_1, ..., \omega_s$ is an $s$-dimensional shift-invariant subspace.

In view of the above, it is now clear that simply belonging to a shift-invariant subspace does not guarantee that a signal $x$ can be reproduced
by a {\em predictive} filter of small $\ell_2$-norm.
For instance, given $r\in \C,\,|r|>1$, consider signals from the parametric family
\[
\X_r=\{x\in\C(\Z):\,x_\tau=\beta r^\tau, \beta \in \C\}.
\]
Here $\X_r$ is a one-dimensional shift-invariant subspace of $\C(\Z)$---solution set of the equation $(1-r\Delta) x= 0$. Clearly, for $x\in \X_r$
$x_t$ cannot be estimated consistently using noisy observations on the left of $t$ (cf.~\cite{shiryaev1997sequential}), and we cannot expect a ``good'' predictive filter to exist for all $x\in \X_r$.

The above example is representative of the difficulties arising when  predicting signals from shift-invariant subspaces of $\C(\Z)$: the characteristic polynomial of the associated difference equation is unstable---its root $z=1/r$ lies {\em inside the (open) unit disk}. Therefore, to be able  to build good ``left'' predictive filters, we need to reduce the class of signals to solutions of equations \rf{eq:shift-inv-diff_eq} with {\em stable polynomials,} with all roots lying outside the (open) unit disk---decaying exponents, harmonic oscillations, and their products. Note that if we are interested in estimating $x_t$ using only observations {\em on the right of $t$}, similar difficulties will arise when $x$ is a solution of a homogeneous linear difference equation with {\em roots outside the closed unit disc}---this situation is completely similar to the above, up to the inversion of the time axis.

\subsubsection{Adaptive prediction of generalized harmonic oscillations}
\label{sec:harmonic-pred}
The above discussion motivates our interest {in} a special family of shift-invariant subspaces which allow for constructing good ``left'' and ``right'' prediction filters---that of sets of solutions to linear homogeneous difference  equations~\rf{eq:shift-inv-diff_eq} with all roots $z_k$  {\em on the unit circle,} i.e.,~$z_k=e^{-i\omega_k}$ with real $\omega_k\in [0,2\pi[$, $k=1,...,s$. In other words,  we are interested in the class of solutions to equation \rf{eq:shift-inv-diff_eq} with $p(z)=\prod_{k=1}^s (1-e^{i\omega_k}z)$ comprised of signals of the form
\[
x_t = \sum_{k=1}^r q_k(t) e^{i \omega_k t}
\]
where $\omega_1, ..., \omega_r \in [0, 2\pi[$ are distinct oscillation frequencies and $q_k(\cdot)$, $k=1,...,r$, are (arbitrary) polynomials of degree $m_k-1$, $m_k$ being the multiplicity of the root $z_k=e^{-i\omega_k}$ (i.e.,~$\sum_{k = 1}^r m_k  = s$).
We call such signals \textit{generalized harmonic oscillations}; we denote $\cH_s[\omega]$ the space of such signals with fixed spectrum $\omega\in [0, 2\pi[^s$ and denote $\cH_s$ the set of generalized harmonic oscillations with at most $s$ (unknown) frequencies.

The problem of constructing a predictive filter for signals from $\cH_s[\omega]$ has already been studied in~\cite{jn-2014}, where the authors proved (cf. \cite[Lemma 6.1]{jn-2014}) that for any $s\geq 1$, vector of frequencies $\omega_1,...,\omega_s$, and $m$ large enough there is $\phi^o \in \C_m(\Z)$ such that
$x=\phi^o*\Delta^mx$ and
\be
\|\phi^o\|_2\leq Cs^{3/2}\sqrt{\log[s+1]\over m}.
\ee{eq:oldl2b}
Here we utilize an improved version of that result.
\begin{prop}
\label{th:sines}
Let $s\geq 1$ and $\omega\in [0,2\pi[^s$. Then for any $m\geq cs^2\log s$ there is a filter $\phi^o\in \C_m(\Z)$ which only depend on $\omega$ such that $x=\phi^o*\Delta^mx$ for all $x\in \cH_s[\omega]$
and
\be
\|\phi^o\|_2 \le Cs\sqrt{\log m\over m}.
\ee{eq:sines-rho}
\end{prop}
\noindent Let now $\cH_{m,n}(s,\varkappa)$ be the set of signals $x\in\C(\Z)$ (locally) close to $\cH_s$ in~$\ell_2$-norm, i.e., which can be decomposed (cf. Assumption~\rf{ass:subspace}) as
\[
x=x{^{\cH}}+\varepsilon
\]
where $x^{\cH}\in\cH_s$ and
\[
\left\|\Delta^{-\tau} \varepsilon\right\|_{n,2} \leq \varkappa \sigma, \quad |\tau|\leq m.
\]
Equipped with the bound of Proposition \ref{th:sines}, we can now derive risk bounds for adaptive predictive estimates on $\cH_{m,n}(s,\varkappa)$.
%For the sake of conciseness, we only formulate these bounds for constrained recovery.
Specifically, following  the proof of Propositions \ref{pr:l2sloss} and \ref{pr:linfsloss} we obtain the following corollaries of the oracle inequalities of Proposition \ref{prop:l2all}.
\begin{prop}
\label{pr:l2sloss-p}
Let $s, m,n \in \Z_+$, $m\geq cs^2\log s$ with large enough $c$, and let $\varkappa\ge 0$.
\item[(i)] Let $\bar \vrho=Cs^2\log m$ with $C$ large enough, and let $\wh x_{\con}=\wh\vphi_{\con}*\Delta^my$ where $\wh\vphi_{\con}$ is an optimal solution  to~\eqref{opt:l2con-pred}; let also $D_n=\{-n+m,...,n+m\}$.
Then for any $\alpha\in ]0,1/2]$
\[
\risk_{D_n,2,\alpha}(\wh x_{\con}|\cH_{m,n}(s,\varkappa))\leq C'\chi^{\alpha}_{m,n}(\sigma,s,\varkappa)
\]
where
\bse
\chi^{\alpha}_{m,n}(\sigma,s,\varkappa)=
\sigma s^2\log [m]\big(\kappa_{m,n}\sqrt{\log[1/\alpha]}+\varkappa\big)
+\sigma s(\kappa_{m,n}+1)\sqrt{\log[ m]\log[(m+n)/\alpha]}.
%\sigma s\sqrt{\log m}\left(s\sqrt{\log m}\big(\kappa_{m,n}\sqrt{\log[1/\alpha]}+\varkappa\big)+(\kappa_{m,n}+1)\sqrt{\log[(m+n)/\alpha]}\right).
\ese
\item[(ii)]
Let $\lambda = \Rem_1^{1/2}$ with $\Rem_1^{\vphantom{1/2}}$ as  in \rf{eq:remainder-pen}, and let $\wh x_{\pen}=\wh\vphi_{\pen}*\Delta^my$ where $\wh\vphi_{\pen}$ is an optimal solution  to~\eqref{opt:l2pen-pred++}.
Then for any $\alpha\in ]0,1/2]$
\[
\risk_{D_n,2,\alpha}(\wh x_{\pen}|\cH_{m,n}(s,\varkappa))\leq C\widetilde\chi^{\alpha}_{m,n}(\sigma,s,\varkappa)
\]
where
\bse
\widetilde\chi^{\alpha}_{m,n}(\sigma,s,\varkappa)&=\sigma s^2\log [m]\left((\kappa_{m,n}+1)\sqrt{\log[(m+n)/\alpha]}+\varkappa\right).
\ese
\end{prop}\noindent
Next, in order to state the result describing pointwise risks of the proposed estimate we need to replace the class $\cH_{m,n}(s,\varkappa)$ with the class of signals which are (locally) ``uniformly'' close to $\cH_s$. Namely, let ${\overline \cH_{m,n}(s,\varkappa)}$ be the set of signals $x\in \C(\Z)$ which can be decomposed (cf.~Assumption~\ref{ass:subspace_infty}) as
\[
x=x^{\cH}+\varepsilon
\]
with  $x{^{\cH}} \in\cH_s$ and
\[
|\varepsilon_\tau| \leq {\varkappa \sigma\over \sqrt{2n+1}}, \quad |\tau|\leq n+m.
\]
\begin{prop}
\label{pr:linfsloss-p}
Let $s, m,n \in \Z_+$, $m\geq cs^2\log s$ with large enough $c$, $n\geq m/2$, and let $\varkappa\ge 0$. We set  $D_{n,m}=\{-n+2m,...,n+m\}$.
\item[(i)] Let $\wh x_{\con}=\wh\vphi_{\con}*\Delta^my$ where $\wh\vphi_{\con}$ is an optimal solution  to~\eqref{opt:l2con-pred} where $\bar \vrho=Cs^2\log m$ with $C$ large enough.
Then for any $\alpha\in ]0,1/2]$
\[
\risk_{{D_{n,m}},\alpha}(\wh x_{\con}|{\overline \cH_{m,n}(s,\varkappa)})\leq C'\nu^\alpha_{m,n}(\sigma,s,\varkappa)
\]
where
\bse \nu^\alpha_{m,n}(\sigma,s,\varkappa)&=&
{s\sqrt{\log m\over m}}\chi^{\alpha}_{m,n}(\sigma,s,\varkappa)
+{\sigma s^3(\log m)^{3/2}\over \sqrt{m}}\left(
\varkappa+\log[1/\alpha]\right)\\
&\leq &C''{\sigma s^3(\log m)^{3/2}\over \sqrt{m}}\left(
\varkappa+\log[1/\alpha]\right).
\ese
\item[(ii)]
Let $\wh x_{\pen}=\wh\vphi_{\pen}*\Delta^my$ where $\wh\vphi_{\pen}$ is an optimal solution  to~\eqref{opt:l2pen-pred++} with $\lambda = \Rem_1^{1/2}$, $\Rem_1^{\vphantom{1/2}}$ being defined in \rf{eq:remainder-pen}.
Then for any $\alpha\in ]0,1/2]$
\[\risk_{{D_{n,m}},\alpha}(\wh x_{\pen}|{\overline \cH_{m,n}(s,\varkappa)})\leq C\widetilde\nu^\alpha_{m,n}(\sigma,s,\varkappa)
\]
where
\bse
\widetilde\nu^\alpha_{m,n}(\sigma,s,\varkappa)&=&
s{\sqrt{\log m\over m}}\widetilde\chi^{\alpha}_{m,n}(\sigma,s,\varkappa)+{\sigma s^3(\log m)^{3/2}\over \sqrt{m}}\left(
\varkappa+\log[1/\alpha]\right)\\\\
&\leq& C'{\sigma s^3(\log m)^{3/2}\over \sqrt{m}}\left(
\varkappa+\log[(m+n)/\alpha]\right).
\ese
\end{prop}

\subsection{Harmonic oscillation denoising}
\label{sec:sines}

To illustrate the results of the previous section, let us consider the problem of recovery of generalized harmonic oscillations. Specifically,  given observations $y_\tau=x_\tau +\sigma \zeta_\tau$, $|\tau|\leq L\in \Z_+$ we  are to estimate the signal $x\in \cH_s$.
We measure the statistical performance of the adaptive estimate $\wh x$ by the maximal over $\cH_s$ integral $\alpha$-risk
\[
\risk_{D_L,2,\alpha}(\wh x|\cH_s)=\inf \left\{r:\sup_{x\in \cH_s}\Prob\left\{\|\wh x-x\|_{L,2}\geq r\right\}\leq \alpha\right\}
\]
on the entire observation domain $D_L=\{-L,...,L\}$.\par
Note that if the frequencies were known, the ordinary least-squares estimate would attain the risk $\cO\left(
\sigma \sqrt{s}\right)$ (up to a logarithmic factor in $\alpha$). When the frequencies are unknown, the lower bound (see, e.g., \cite[Theorem 2]{recht2}) states that
\begin{equation}
\label{eq:sines-lower}
\risk_{D_L,2,\frac{1}{2}}(\wh x|{\cH_s}) \ge c\sigma \sqrt{s \log L}.
\end{equation}
In the case where all frequencies are different, this bound is attained asymptotically by the maximum likelihood estimate \cite{tufts1982estimation,stoica1989music}. However, implementing that estimate involves computing maximal likelihood estimate of $\omega$---a global minimizer in the
optimization problem
\[
\min_{\alpha\in \C^s,\,\omega\in \R^s}\left(\sum_{|\tau|\leq L}\left|y_\tau-\sum_{k=1}^s \alpha_k e^{i \omega_k \tau}\right|^2\right)^{1/2}
\]
 and becomes numerically challenging already for very moderate values of $s$.
Moreover, the lower bound ~\eqref{eq:sines-lower} is in fact attained by the \textit{Atomic Soft Thresholding} (AST) estimate \cite{recht1,recht2}---which can be implemented efficiently---but only under the assumption that the frequencies $\{\omega_1, ..., \omega_s\}$ are well separated---precisely, when the {\em minimal frequency separation} in the wrap-around distance
%matching upper bound for the AST is only proved (cf. \cite[Theorem 1]{recht2})
\begin{equation}
\label{def:min-sep}
\delta_{\min}  := \min_{1 \le j \ne k \le s}{\min \{ |\omega_j - \omega_{k}|, 2\pi - |\omega_j - \omega_{k}|\}}
\end{equation}
satisfies $\delta_{\min}> {2\pi\over 2L+1}$ (cf. \cite[Theorem 1]{recht2}).
To the best of our knowledge, the question whether there exists an efficiently implementable estimate matching the lower bound~\eqref{eq:sines-lower} in the general case is open.

A new approach to the problem was suggested in~\cite{harchaoui2015adaptive} where a uniform-fit adaptive estimate was used for estimation and prediction of (generalized) harmonic oscillations. That approach, using the bound \rf{eq:oldl2b} along with the estimate for the risk of the uniform-fit recovery, resulted in the final risk bound $\cO\left(\sigma s^3\log[s]\log [L/\alpha]\right)$.

Using the results in the preceding section we can now build an improved adaptive estimate. Here we assume that {the number~$s$ of frequencies (counting with their multiplicities)} is known in advance,  and  utilize constrained recoveries \rf{opt:l2con} and \rf{opt:l2con-pred} with the parameter $\bar\vrho$ selected using this information;\footnote{It is worth mentioning that the AST estimate does not require the a priori knowledge of $s$;  we can also get rid of this hypothesis when using the procedure which is adaptive to the unknown value of $s$, at the expense of an additional logarithmic factor.} note that~$s$ is precisely the dimension of the shift-invariant subspace to which~$x$ belongs, cf.~Proposition~\ref{th:shift_invariant}. Let us consider the following procedure.
  \begin{quotation}\noindent
Choose $K \le L$, and divide the observation interval~$D_L$ into the central segment~$D_K=\{-K,...,K\}$
and left and right segments $D_- = \{-L,...,-K-1\}$ and $D_+ =\{K+1,...,L\}$. In what follows we assume that $L$ and $K$ are even and put $k=(L-K)/2$. Then we act as follows.
\begin{itemize}
\item Using the data $y_\tau,\,|\tau|\leq L$ we compute an optimal solution $\wh\vphi \in \C_{L-K}(\Z)$ to the optimization problem \rf{opt:l2con} with $m=L-K$, $n=K$, and  $\bar \vrho=4s$; for  $t\in D_n$ we compute the interpolating (two-sided) estimate $\wh x_t =[\wh\vphi *y]_t$.
\item We set $m=\lfloor (L+n)/2\rfloor$, $n=k$, $\bar\vrho=\bar\vrho^{+} := 2C^2s^2\log L$ where $C$ is as in the bound \rf{eq:sines-rho} of Proposition \ref{th:sines} and compute an optimal solution $\wh\vphi^+ \in \C_{m}(\Z)$ to the optimization problem \rf{opt:l2con-pred};
for~$t\in D_+$ we compute the left (one-sided)  prediction $\wh x_t =[\wh\vphi^+*\Delta^my]_t$.
\item We set $m=\lfloor (L+n)/n\rfloor$, $n=k$,  $\bar\vrho=\bar\vrho^{+}$ and compute an optimal solution $\wh\vphi^- \in \C_{m}(\Z)$ to the ``right'' analog of \rf{opt:l2con-pred};\footnote{In
    the corresponding ``right'' optimization problem the ``left prediction'' $\vphi *\Delta^my$ is replaced with the ``right prediction'' $\vphi *\Delta^{-m}y$. Therefore, the objective to be minimized in this case is
    $\|\Delta^{m}(y-\vphi *\Delta^{-m} {y})\|_{n,2}$.}  for~$t\in D_-$ we compute the right (one-sided) prediction $\wh x_t =[\wh\vphi^-*\Delta^{-m}y]_t$ .
\end{itemize}
We select $K$ to minimize the ``total'' risk bound of the adaptive recovery over $D_L$.
  \end{quotation}
We have the following corollary of the Propositions \ref{pr:l2sloss} and \ref{pr:l2sloss-p} in the present setting.
\begin{prop}
\label{pr:sines}
Suppose that $L\geq c s^2\log s$ with large enough $c>0$. Then, in the situation of this section, for any $\alpha\in]0,1/2]$
\be
\risk_{D_L,2,\alpha}(\wh x|\cH_s) \leq C \sigma s^{3/2}\log [L/\alpha] .
\ee{eq:sines-composite}
\end{prop}
\paragraph{Remarks.}
The risk bound \rf{eq:sines-composite}, while significantly improved in terms of dependence {on} $s$ over the corresponding bound of \cite{harchaoui2015adaptive}, contains an extra factor $O(s\sqrt{\log L})$ {when} compared to the lower bound \rf{eq:sines-lower}.
It is unclear to us whether this factor can be reduced for an efficiently computable estimate.
\par
It may be worth mentioning  that when the frequency separation assumption holds, i.e.,~when $\delta_{\min}>{2\pi\over 2L+1}$ where the separation $\delta_{\min}$ is defined in \rf{def:min-sep}, the above estimation procedure can be simplified: one can ``remove'' the central segment in the above construction only using left and right adaptive predictive estimates on two half-domains. The ``total'' $(1-\alpha)$-reliable  $\ell_2$-loss of the ``simplified'' adaptive recovery is then \[
\cO\left(\sigma \sqrt{s^2\log [1/\alpha] + s\log[L/\alpha]}\right).
\]
The latter bound is a simple corollary of the oracle inequalities of Proposition \ref{prop:l2all} and the following statement.
 \begin{lemma}
\label{lem:sepsines}
Let $m\in \Z_+$, $\nu > 1$, and let $\cH_s[\omega]$ be the set of harmonic oscillations $x$ with the minimal frequency separation satisfying
\begin{equation}
\label{eq:min-sep}
\delta_{\min} \ge \frac{2\pi \nu}{2m+1}.
\end{equation}
Then there exists a filter $\phi^o \in \C_{m}(\Z)$ satisfying $x = \phi^o * \Delta^mx$ for all $x\in \cH_s[\omega]$ and such that
\begin{equation*}
\|\phi^o\|_2 \le \sqrt{\frac{Qs}{2m+1}}, \quad \text{where} \quad Q = \frac{\nu+1}{\nu-1}.
\end{equation*}
In particular, whenever
$
\delta_{\min} \ge {{4\pi}\over {2m+1}},
$
one has
\[
\|\phi^o\|_2 \le \sqrt{\frac{3s}{2m+1}}.
\]
\end{lemma}

\section*{Acknowledgements}
 Dmitrii Ostrovskii was supported by ERCIM Alain Bensoussan Scholarship while finalizing this project.
Zaid Harchaoui received support from  NSF CCF 1740551.
Research of Anatoli Juditsky and Arkadi Nemirovski was supported by MIAI \@ Grenoble Alpes (ANR-19-P3IA-0003).

\appendix
\section{Preliminaries}
First, let us present some additional notation and technical tools to be used in the proofs.

\subsection{Additional notation}\label{sec:addnot}
In what follows, $\Re(z)$ and $\Im(z)$ denote, correspondingly, the real and imaginary parts of  $z \in \C$, and $\overline{z}$ denotes the complex conjugate of $z$. For a matrix $A$ with complex entries, $\overline{A}$ stands for the conjugation of $A$ (without transposition), $A^\T$ for the transpose of $A$, and $A^\H$ for its Hermitian conjugate. We denote $A^{-1}$ the inverse of $A$ when it exists.
$\Tr(A)$ denotes the trace of a matrix $A$ and $\det A$ its determinant; $\|A\|_\F$ is the Frobenius norm of $A$,  $\|A\|_*$ is the operator norm, and $\|A\|$ is the nuclear norm.
We also denote $\lmax(A)$ and $\lmin(A)$ the maximal and minimal eigenvalues of a Hermitian matrix $A$.
For  $a \in \C^{n}$ we denote $\Diag(a)$ the $n\times n$ diagonal matrix with diagonal entries $a_i$.
%We denote $I$ the identity matrix, sometimes with a subscript indicating its size.
We use notation $\|x\|^*_{n,p}$ for the $\ell_p$-norm of the DFT of $x$ so that
\[
\|x\|^*_{n,p}=\|F_n[x]\|_p=\left(\sum_{k=1}^{2n+1} \big|\big(F_n[x]\big)_k\big|^p\right)^{1/p}
\]
with the standard interpretation of $\|\cdot\|^*_{n,\infty}$.

In what follows, we associate linear maps $\C_{n}(\Z) \to \C_{n'}(\Z)$ with matrices in~$\C^{(2n+1)\times (2n'+1)}$.

%matrices with complex elements.

\paragraph{Convolution matrices.} We use the following  matrix-vector representations of discrete convolution.
\begin{itemize}
\item
Given $y \in \C(\Z)$, we associate with it an $(2n+1)\times(2m+1)$ matrix
\be
T(y)=
%\underbrace{
\left[\begin{array}{clclc}
y_{-n+m}& \cdots & y_{-n}&\cdots &y_{-n-m}\\
\vdots &\cdots & \vdots &\cdots &\vdots \\
y_m &  \cdots & y_{0}& \cdots &y_{-m}\\
%y_{m+1} & \cdots & y_{1}&\cdots &y_{1-m}\\
\vdots &\cdots & \vdots &\cdots &\vdots \\
y_{n+m}& \cdots & y_{n}&\cdots &y_{n-m}
\end{array}\right],
\ee{s(psi)}
such that $[\varphi*y]_{-n}^n = T(y) [\varphi]_{-m}^m$ for $\varphi \in \C_m(\Z)$. Its squared Frobenius norm satisfies
\be
\left\|T(y)\right\|_\F^2 = \sum\limits_{|\tau|\leq m} \|\Delta^\tau y\|^2_{n,2}.
\ee{tr21}
\item
Given $\varphi \in \C_m(\Z)$, consider a $(2n+1)\times(2m+2n+1)$ matrix
\be
M(\varphi)=
\left[\begin{array}{cllllllc}
\varphi_{m}&\cdots &\cdots  &\varphi_{-m}&0&\cdots &\cdots &0\\
0&\varphi_{m}&\cdots  &\cdots &\varphi_{-m}&0&\cdots &0\\
\vdots &\ddots &\ddots &\cdots &\cdots &\ddots &\cdots &\vdots \\
\vdots &\cdots &\ddots &\ddots &\cdots &\cdots &\ddots &\vdots \\
0 &\cdots &\cdots &0&\varphi_{m}&\cdots &\cdots &\varphi_{-m}
\end{array}\right],
\ee{M(phi)}
such that for $y \in \C(\Z)$ one has $[\varphi*y]_{-n}^n=M(\varphi) [y]_{-m-n}^{m+n}$, and
\be
\left\|M(\varphi)\right\|_\F^2= (2n+1)  \|\varphi\|^2_{m,2}.
\ee{tr22}
%For instance, when $\varphi$ is uniformly bounded, $\|\varphi\|_\infty \leq b$, one has $\|\varphi\|^2_{m,2}\leq (m+1)b^2$ and
%\begin{equation*}
%\left\|M(\varphi)\right\|_F^2 \leq (n+1)(m+1)b^2.
%\end{equation*}
\item
Given $\varphi \in \C_m(\Z)$, consider the following circulant matrix of size $2m+2n+1$:
\be
C(\varphi) =
\left[\begin{array}{ccccccccccccc}
\varphi_{0} &\cdots & \cdots & \varphi_{-m} & 0 & \cdots & \cdots & \cdots &  0  & \vphi_m & \cdots & \cdots & \vphi_1 \\
\vphi_1 & \varphi_{0} &\cdots & \cdots & \varphi_{-m} & 0 & \cdots & \cdots & \cdots &  0  & \vphi_m & \cdots & \vphi_2 \\
\cdots & \cdots & \ddots & \cdots & \cdots & \ddots & \ddots & \cdots & \cdots & \cdots & \cdots & \cdots & \cdots \\
\cdots & \cdots & \cdots & \ddots & \cdots & \cdots & \ddots & \ddots & \cdots & \cdots & \cdots & \cdots & \cdots \\
\cdots & \cdots & \cdots & \cdots & \ddots & \cdots & \cdots & \ddots & \ddots & \cdots & \cdots & \cdots & \cdots \\
\cdots & \cdots & \cdots & \cdots & \cdots & \ddots & \cdots & \cdots & \ddots & \ddots & \cdots & \cdots & \cdots \\
0 & \cdots  & 0 & \varphi_{m} &\cdots & \cdots & \varphi_{0} & \cdots & \cdots & \vphi_{-m} & 0 & \cdots & 0 \\
\cdots & \cdots & \cdots & \ddots & \ddots & \cdots & \cdots & \ddots & \cdots & \cdots & \cdots & \cdots & \cdots \\
\cdots & \cdots & \cdots & \cdots & \ddots & \ddots & \cdots & \cdots & \ddots & \cdots & \cdots & \cdots & \cdots \\
\cdots & \cdots & \cdots & \cdots & \cdots & \ddots & \ddots & \cdots & \cdots & \ddots & \cdots & \cdots & \cdots \\
\cdots & \cdots & \cdots & \cdots & \cdots & \cdots & \ddots & \ddots & \cdots & \cdots & \ddots & \cdots & \cdots \\
\cdots & \cdots & \cdots & \cdots & \cdots & \cdots & \cdots & \ddots & \ddots & \cdots & \cdots & \ddots & \cdots \\
\varphi_{-1} &\cdots & \cdots & \varphi_{-m} & 0 & \cdots & \cdots  & \cdots &  0  & \vphi_m & \cdots & \cdots & \vphi_0
\end{array}\right].
%\left[\begin{array}{ccccccccccccccc}
%\varphi_{0} &\cdots & \cdots & \cdots & \varphi_{-m} & 0 & \cdots & \cdots  & \cdots & \cdots &  0  & \vphi_m & \cdots & \cdots & \vphi_1 \\
%\vphi_1 & \varphi_{0} &\cdots & \cdots & \cdots  & \varphi_{-m} & 0 & \cdots & \cdots & \cdots & \cdots &  0  & \vphi_m & \cdots & \vphi_2 \\
%\cdots & \cdots &\ddots & \cdots & \cdots  & \cdots & \ddots & \cdots & \cdots & \cdots & \cdots &  \cdots  & \cdots & \cdots & \cdots \\
%\cdots & \cdots &\cdots & \ddots & \cdots  & \cdots & \cdots & \ddots & \cdots & \cdots & \cdots &  \cdots  & \cdots & \cdots & \cdots \\
%\cdots & \cdots &\cdots & \cdots & \ddots  & \cdots & \cdots & \cdots & \ddots & \cdots & \cdots &  \cdots  & \cdots & \cdots & \cdots \\
%\cdots & \cdots &\cdots & \cdots & \cdots  & \ddots & \cdots & \cdots & \cdots & \ddots & \cdots &  \cdots  & \cdots & \cdots & \cdots \\
%\cdots & \cdots &\cdots & \cdots & \cdots  & \cdots & \ddots & \cdots & \cdots & \cdots & \ddots &  \cdots  & \cdots & \cdots & \cdots \\
%0 & \cdots  & 0 & \varphi_{m} &\cdots & \cdots & \cdots & \varphi_{0} & \cdots & \cdots & \cdots & \vphi_{-m} & 0 & \cdots & 0 \\
%\cdots & \cdots &\cdots & \cdots & \cdots  & \cdots & \cdots & \cdots & \ddots & \cdots & \cdots &  \cdots  & \ddots & \cdots & \cdots \\
%\cdots & \cdots &\cdots & \cdots & \cdots  & \cdots & \cdots & \cdots & \cdots & \ddots & \cdots &  \cdots  & \cdots & \ddots & \cdots \\
%\varphi_{-1} &\cdots & \cdots & \varphi_{-m} & 0 & \cdots & \cdots  & \cdots & \cdots &  0  & \vphi_m & \cdots & \cdots & \cdots & \vphi_0
%\end{array}\right].
\ee{C(phi)}
Note that~$C(\varphi) [y]_{-m-n}^{m+n}$ is the circular convolution of $[y]_{-m-n}^{m+n}$ and the zero-padded filter
\[
\tilde\vphi := [\vphi]_{-m-n}^{m+n} = [0; ...; \vphi_{-m}; ...; \vphi_m; 0; ...; 0],
\]
that is, convolution of the periodic extensions of $[y]_{-m-n}^{m+n}$ and $\tilde \vphi$ evaluated on $\{-m-n,...,m+n\}$.
Hence, by the diagonalization property of the DFT operator one has
\begin{equation}\label{eq:circ_diag}
C(\varphi) = \sqrt{2m+2n+1}  F^\H_{m+n} \text{diag}(F_{m+n}\tilde\vphi) F_{m+n}
\end{equation}
where with some notational abuse we denote $F_{n}$ the matrix of DFT with the entries
\[
[F_n]_{kj}={1\over \sqrt{2n+1}}\exp\left({2\pi i(k-n)j\over 2n+1}\right),\quad 1\leq k,j\leq 2n+1.\]
Besides this, note that
\begin{equation*}
\left\|C(\varphi)\right\|_\F^2 = (2m+2n+1) \|\varphi\|^2_{m,2}.
\end{equation*}
%This matrix is useful since
\end{itemize}

\paragraph{Reformulation of approximate shift-invariance}
The following reformulation of Assumption~\ref{ass:subspace} will be convenient for our purposes.
\begin{quotation}\noindent
{\em There exists an $s$-dimensional vector subspace $\S_{n}$ of $\C^{2n+1}$
%consisting of elements of $\mathcal{S}$ restricted on any $\C\big(\Z_{-\tau}^{-\tau+n}\big)$, $0\le\tau \le m$,
and an idempotent Hermitian $(2n+1)\times (2n+1)$ matrix~$\Pi_{\S_{n}}$ of rank {$s$---}projector on $\S_n$---such that
\be
\big\|\left(I_{2n+1}-\Pi_{\mathcal{S}_{n}}\right) \left[\Delta^\tau x\right]_{-n}^{n}\big\|_2 \, \Big[={\left\|\Delta^{\tau} \varepsilon\right\|_{n,2}
}\Big] \leq \sigma\varkappa, \quad |\tau|\leq m
\ee{subspace1}
where $I_{2n+1}$ is the $(2n+1)\times (2n+1)$ identity matrix.
}
\end{quotation}
\subsection{Technical tools}

\paragraph{Deviation bounds for quadratic forms.}
Let $\zeta \sim \CN(0,I_n)$ be a standard complex Gaussian vector, meaning that $\zeta = \xi_1 + i\xi_2$ where $\xi_{1}$  and $\xi_{2}$ are two independent draws from $\N(0,I_n)$. We use simple facts listed below.
\begin{itemize}
\item Due to the unitarity of the DFT, if $\zeta_{-n}^n \sim \CN(0,I_{2n+1})$ we also have $F_n[\zeta]\sim \CN(0,I_{2n+1})$.
\item
 We use a simple bound
\begin{align}\label{eq:complex_gaussian_max}
\prob\left\{\|\zeta\|_{n,\infty} \le \sqrt{2\log n + 2u}\right\} \ge 1-e^{-u}
\end{align}
which can be verified directly using that~$|\zeta_1|^2_2 \sim \chi_2^2$.
\item
The following deviation bounds for $\|\zeta\|^2_2 \sim \chi^2_{2n}$ are due to \cite[Lemma 1]{lama2000}:
\begin{equation}\label{eq:chi_bound_low}
\begin{aligned}
&\prob\left\{\frac{\|\zeta\|^2_2}{2} \le n + \sqrt{2nu} + u \right\} \ge 1-e^{-u},\\
&\prob\left\{\frac{\|\zeta\|^2_2}{2} \ge n - \sqrt{2nu} \right\} \ge 1-e^{-u}.
\end{aligned}
\end{equation}
By simple algebra we obtain an upper bound for the norm:
\begin{align}\label{eq:chi_bound}
&\prob\left\{\|\zeta\|_2 \le \sqrt{2n} + \sqrt{2u} \right\} \ge 1-e^{-u}.
\end{align}
\item
Further, let $K$ be an $n \times n$ Hermitian matrix with the vector of eigenvalues $\lambda = [\lambda_1;\, ...;\, \lambda_n]$.
Then the real-valued quadratic form $\zeta^\H K \zeta$ has the same distribution as $\xi^T B\xi$, where $\xi = [\xi_1; \xi_2] \sim \N(0,I_{2n})$, and $B$ is a real $2n \times 2n$ symmetric matrix with the vector of eigenvalues $[\lambda; \lambda]$.
We have $\text{Tr}(B) = 2\text{Tr}(K)$, $\|B\|^2_\F = 2\|K\|_\F^2$ and $\|B\| = \|K\| \le \|K\|_\F$.
Invoking again~\cite[Lemma 1]{lama2000} (a close inspection of the proof shows that the assumption of positive semidefiniteness can be relaxed), we have
\begin{align}\label{trq2}
\prob\left\{\frac{\zeta^\H K \zeta}{2} \leq \Tr(K) + (u + \sqrt{2u})\|K\|_\F\right\} \geq 1-e^{-u}.
\end{align}
Further, when $K$ is positive semidefinite, we have $\|K\|_\F \le \text{Tr}(K)$, whence
\begin{align}\label{trq1}
\prob\left\{\frac{\zeta^\H K \zeta}{2} \leq \Tr(K)(1+\sqrt{u})^2\right\} \geq 1-e^{-u}.
\end{align}
\end{itemize}
The following lemma, interesting in its own right, controls the inflation of the $\ell_1$-norm of the DFT of a zero-padded signal.
\begin{lemma}
\label{lemma:zero_padding}
Let $u \in \C_m(\Z)$ one has
%for $\kappa_{m,n} := \sqrt{\frac{n+1}{m+1}}$,
\begin{equation*}
\|u\|^*_{m+n,1} \le \|u\|^*_{m,1} (1 + \kappa^2_{m,n})^{1/2} [\log(m+n+1)+3].
\end{equation*}
\end{lemma}
\paragraph{Proof.} It suffices to show that the bound
\[\|u\|^*_{m+n,1} \le (1 + \kappa^2_{m,n})^{1/2} [\log(m+n+1)+3]
\]
holds for all $u\in \C_m(\Z)$ such that $\|u\|^*_{m,1} \le 1$. We assume that $n \ge 1$, the lemma statement being trivial otherwise.

First of all, function $\|u\|^*_{m+n,1}$ is convex so its maximum over the set $u\in \C_m(\Z)$, $\|u\|^*_{m,1} \le 1$, is attained at an extreme point $u^j$ of the set given by
$F_m [u^j] = e^{i \theta} {e}^j$ where ${e}^j$ is the $j$-th canonic basis vector and $\theta \in [0,2\pi]$.
Note that
\[
u^j_\tau = \frac{1}{\sqrt{2m+1}} \exp\left(i\left[\theta+\frac{2\pi \tau j}{2m+1}\right]\right),
\]
thus, for~$\gamma_{m,n} := \sqrt{(2m+2n+1)(2m+1)}$ we obtain
\begin{align*}
\big\|u^j\big\|^*_{m+n,1}
&= \frac{1}{\gamma_{m,n}} \sum_{k=1}^{2(m+n)+1} \left| \sum_{|\tau| \leq m} \exp\left(2\pi i \tau\left[\frac{j}{2m+1} - \frac{k}{2m+2n+1} \right]\right) \right| \\
&= \frac{1}{\gamma_{m,n}} \sum_{k=1}^{2(m+n)+1} \left| {\DirKer}_m\left(\omega_{jk}\right)\right|,
\end{align*}
where
\[
\omega_{jk} := 2\pi \left[\frac{j}{2m+1} - \frac{k}{2m+2n+1} \right]
\]
and
${\DirKer}_m(\cdot)$ is the Dirichlet kernel of order $m$:
\begin{equation*}
{\DirKer}_m(\omega) := \left\{
\begin{aligned}
&\frac{\sin \left({(2m+1)\omega}/{2}\right)}{\sin\left({\omega}/{2}\right)}, & \quad \omega \ne 2\pi l,\\
&2m+1, & \quad \omega = 2\pi l.
\end{aligned}\right.
\end{equation*}
Hence,
\begin{align}\label{eq:dirichlet_sum}
\gamma_{m,n} \|u^j\|^*_{m+n,1} \le \max_{\theta \in [0,2\pi]} \left\{  \Sigma_{m,n}(\theta) := \sum_{k=1}^{2(m+n)+1} \left| {\DirKer}_m\left(\frac{2\pi k}{2m+2n+1}  + \theta \right)\right| \right\}.
\end{align}
%Note that
%\begin{equation}\label{dirichlet_grid}
%\begin{aligned}
%S(0) \le  m+1 + \sum_{k=1}^{n+T} \frac{1}{\left| \sin\left( \frac{\pi k}{n+m+1} \right)\right|}
%&\le m+1 + 2\, \sum_{k=1}^{\left\lceil \frac{n+T}{2} \right\rceil} \frac{1}{\left| \sin\left( \frac{\pi k}{n+m+1} \right)\right|} \\
%&\le m+1 + (n+m+1) \sum_{k=1}^{\left\lceil \frac{n+T}{2} \right\rceil} \frac{1}{k} \\
%&\le (n+m+1)\left[\log(n+m+1) + 1\right],
%\end{aligned}
%\end{equation}
%Note that $|{\DirKer}_m(x)|$ on the unit circle is bounded from above by a positive function~\dimacorr{}{draw a pic!}
%\begin{align*}
%U_m(x) = \left\{
%\begin{array}{ll}
%\frac{\pi}{2 \min (x, \pi-x)}, & \quad x \in (0,\pi),\\
%2m+1, & \quad x = 0.
%\end{array}
%\right.
%\end{align*}
For any $\theta \in [0,2\pi]$, the summation in \eqref{eq:dirichlet_sum} is over the $\theta$-shifted regular $(2m+2n+1)$-grid on the unit circle.
The contribution to the sum $\Sigma_{m,n}(\theta)$ of the two closest to $x = 1$ points of this grid is at most $2(2m+1)$.
Using the bound
\[
{\DirKer}_m(\omega) \le |\sin(\omega/2)|^{-1} \le \frac{\pi}{\min(\omega, 2\pi-\omega)}.
\]
 for the remaining points, and because $f(\omega) = \frac{\pi}{\omega}$ is decreasing on $[\frac{2\pi}{2m+2n+1}, \pi]$ (recall that $n \ge 1$) we arrive at the bound
 \begin{align*}
\Sigma_{m,n}(\theta)
&\le 2\left(2m+1 + \sum_{k=1}^{m+n+1} \frac{2m+2n+1}{2k}\right).
\end{align*}
Now, using the inequality  $H_n \leq \log n + 1$ for the $n$-th harmonic number we arrive at the bound
\bse
\Sigma_{m,n}(\theta) &\le& 2(2m+1) + (2m+2n+1) \left[\log(m+n+1)+1\right] \\&\le& (2m+2n+1) \left[\log(m+n+1)+3\right]
\ese
which implies the lemma. \qed

\section{Proof of Theorems \ref{th:l2con} and \ref{th:l2pen++}}
\label{sec:orc-proofs}
\paragraph{What is ahead.} While it is difficult to describe informally the ideas underlying the proofs of the oracle inequalities,  the ``mechanics'' of the proof of inequality \rf{eq:th1},  for instance, is fairly simple: for any $\vphi^o$ which is feasible to \rf{opt:l2con} one has
\[
\|y-\wh\varphi * y\|_{n,2}\leq \|y-\vphi^o* y\|_{n,2},
\]
and to prove the inequality \rf{eq:th1} all we need to do is to bound tediously all terms of the remainder $\|x-\wh\varphi * y\|_{n,2}-\|x-\varphi^o * y\|_{n,2}$. This may be compared to bounding the $\ell_2$-loss of the Lasso regression estimate. Indeed, let $m = n$ for simplicity, and, given $y\in\C(\Z)$, let  $T(y)$ be the $(2n+1)\times(2n+1)$ ``convolution matrix'' as defined by~\eqref{s(psi)}
%\bse
%T(y)=
%%\underbrace{
%\left[\begin{array}{llll}
%y_{0}&y_{-1}&...&y_{-n}\\
%y_{1}&y_{0}&...&y_{1-n}\\
%...&...&...&...\\
%y_{n}&y_{n-1}&...&y_{0}
%\end{array}\right].
%%}_{2L+2T+1}\;
%\ese{}
such that for $\varphi \in \C_n(\Z)$ one has $[\varphi*y]_{0}^n = T(y) [\varphi]_{-n}^n$.
When denoting $f=F_n[\varphi]$, the optimization problem in \eqref{opt:l2con} can be recast as a ``standard'' $\ell_1$-constrained least-squares problem with respect to $f$:
\be
\min_{f\in \C^{2n+1}} \| y-A_n f\|_{n,2}^2 \;\; \text{s.t.} \;\; \|f\|_1\leq \frac{\bar{\vrho}}{\sqrt{2n+1}}
\ee{eq:f2p1}
where $A_n=T(y)F^{\H}_n$.
Observe that $f^o=F_n [\varphi^o]$ is feasible for \eqref{eq:f2p1}, so that
\[
\|y-A_n\wh{f}\|_{n,2}^2\leq \|y-A_n f^o\|_{n,2}^2,
\]
where $\wh{f}=F_n[\wh{\varphi}]$, and
\[
\begin{aligned}
\|x-A_n\wh{f}\|_{n,2}^2-\|x-A_nf^o\|_{n,2}^2
&\leq 2\sigma \big(\Re\langle\zeta,x-A_nf^o\rangle_n -\Re\langle \zeta,x-A_n\wh{f}\rangle_n\big)\nn
&\leq 2\sigma \big|\langle\zeta,A_n(f^o-\wh{f})\rangle_n\big|\leq 2\sigma \|A_n^\H[\zeta]_{-n}^n\|_\infty\|f^o-\wh{f}\|_1 \nn
&\leq 4\sigma \|A_n^\H[\zeta]_{-n}^n\|_\infty{\bar{\varrho}\over \sqrt{n+1}}.
\end{aligned}
\]
In the ``classical'' situation, where $[\zeta]_{-n}^n$ is independent of $A_n$ (see, e.g., \cite{judnem2000}) one would have
\[
\|A_n^\H [\zeta]_{-n}^n\|_\infty \le c_\alpha \sqrt{\log n}\max_{j}\|[A_n]_j\|_2 \leq c_\alpha \sqrt{n\log n}\max_{i,j}|A_{ij}|
\]
where
%$\|A\|_\infty=\max_{i,j}|A_{ij}|$, and
$c_\alpha$ is a logarithmic in $\alpha^{-1}$ factor. This would rapidly lead to the bound equivalent to \eqref{eq:th1}.
The principal difference with the standard setting which is also the source of the main difficulty in the analysis of the properties of adaptive estimates is that the ``regression matrix'' $A_n$ in the case we are interested in is built of the noisy observations $[y]_{-n}^n$ and thus depends on $[\zeta]_{-n}^n$. In this situation, curbing the cross term is more involved and calls for Assumption~\ref{ass:subspace}.

\subsection{Proof of Theorem~\ref{th:l2con}}

\paragraph{1$^o$.}
Let~$\vphi^o \in \C_m(\Z)$ be any filter satisfying the constraint in~\eqref{opt:l2con}.
Then,
\be
\|x-\wh\varphi*y\|_{n,2}^2
\hspn &\le& \hspn \|(1-\varphi^{o})*y\|_{n,2}^2-{\sigma^2\|\zeta}\|_{n,2}^2-2\sigma\Re \langle \zeta, x-\wh\varphi*y\rangle_n \nn
\hspn &=& \hspn\|x-\varphi^{o}*y\|_{n,2}^2  - 2\underbrace{\sigma\Re\langle \zeta, x-\wh\varphi*y\rangle_n }_{{\delta^{(1)}}} + 2\underbrace{\sigma \Re\langle \zeta, x-\varphi^{o}*y\rangle_n }_{{\delta^{(2)}}}.
\ee{first0}
%\mynote{Dima, is notation $\langle...\rangle_n$ defined somewhere? How do you define it? And what if the elements are "shifted"? You need another notation for that too?}
Let us bound ${\delta^{(1)}}$.
Denote for brevity $I := I_{2n+1}$, and recall that $\Pi_{\cS_n}$ is the projector on $\mathcal{S}_n$ from \eqref{subspace1}. We have the following decomposition:
\begin{equation}
\label{I_pi}
\begin{aligned}
{\delta^{(1)}}
&=\underbrace{\sigma \Re\langle [\zeta]_{-n}^n, \Pi_{\cS_n}[x-\wh\varphi*y]_{-n}^n\rangle}_{{\delta^{(1)}_1}} +\underbrace{\sigma \Re\langle [\zeta]_{-n}^n, (I-\Pi_{\cS_n})[x-\wh\varphi*x]_{-n}^n\rangle}_{{\delta^{(1)}_2}} \nn
&\quad-\underbrace{\sigma^2 \Re\langle [\zeta]_{-n}^n, (I-\Pi_{\cS_n})[\wh\varphi * {\zeta}]_{-n}^n \rangle}_{{\delta^{(1)}_3}}
\end{aligned}
\end{equation}
One can easily bound $\delta^{(1)}_1$ under the premise of the theorem:
\[
\begin{aligned}
\left|{\delta^{(1)}_1} \right|%\leq|(z_{-n}^n)^*\Pi_{\cS_n}[x-\wh\varphi*y]_0^n|
&\leq {\sigma}\big\|\Pi_{\cS_n} {[\zeta]}_{-n}^n\big\|_2  \big\|\Pi_{\cS_n}[x-\wh\varphi*y]_{-n}^n\big\|_{2} \nn
&\leq {\sigma}\big\|\Pi_{\cS_n} {[\zeta]}_{-n}^n\big\|_2  \big\|x-\wh\varphi*y\big\|_{n,2}.
\end{aligned}
\]
Note that $\Pi_{\cS_n} [\zeta]_{-n}^n \sim \CN(0,I_s)$, and by \eqref{eq:chi_bound} we have
\[
\prob\left\{\big\|\Pi_{\cS_n}[\zeta]_{-n}^n\big\|_2 \geq \sqrt{2s}+\sqrt{2u}\right\}\leq e^{-u},
\]
which gives the bound
\be
\prob\left\{\big|{\delta^{(1)}_1}\big|%|(\zeta_{-L}^L)^*\Pi_{\cS_n}(x-\wh\varphi*y)_{-L}^L|
\leq \sigma\big\|x-\wh\varphi*y\big\|_{n,2} \left(\sqrt{2s}+\sqrt{2\log\left[1/\alpha_1\right]}\right) \right\} \ge 1-\alpha_1.
\ee{i11}

\paragraph{2$^o$.} We are to bound the second term of \eqref{I_pi}. To this end, note first that
\bse
{\delta^{(1)}_2} = \sigma \Re \langle [\zeta]_{-n}^n, (I-\Pi_{\cS_n}) [x]_{-n}^n \rangle - \sigma \Re \langle [\zeta]_{-n}^n, (I-\Pi_{\cS_n})[\wh\varphi*x]_{-n}^n\rangle.
\ese
By \eqref{subspace1}, $\left\|(I-\Pi_{\cS_n})[x]_{-n}^n\right\|_{2}\leq \sigma \varkappa$, thus with probability $1-\alpha$,
\be
\left|{\langle [\zeta]_{-n}^n}, (I-\Pi_{\cS_n})[x]_{-n}^n\rangle \right| \le \sigma\varkappa\sqrt{2\log[1/\alpha]}.
\ee{i22-1}
On the other hand, using the notation defined in \eqref{s(psi)}, we have
 $[\wh\varphi*x]_{-n}^n=T(x)[\wh\varphi]_{-m}^m$, so that
\[
{\langle [\zeta]_{-n}^n}, (I-\Pi_{\cS_n})[\wh\varphi*x]_{-n}^n\rangle = {\langle [\zeta]_{-n}^n}, (I-\Pi_{\cS_n})T(x)[\wh\varphi]_{-m}^m\rangle.
\]
%and
%\[
%\left|(\zeta_0^n)^*(I-\Pi_{\cS_n})(\wh\varphi*x)_0^L\right|\leq \|S^*(x)(I-\Pi_{\cS_n})(\zeta_0^n)^*\|_{2}\|\wh\varphi_{-T}^T\|_2.
%\]
Note that $[T(x)]_{\tau} = [\Delta^{\tau} x]^{n}_{-n}$ for the columns of $T(x)$, $|\tau|\leq m$. By \eqref{subspace1}, we have
%\dimacorr{$$(I-\Pi_{\cS_n}) T(x)[\tau+1] = (I-\Pi_{\cS_n})x^{-\tau+n}_{-\tau}=\varepsilon^{\tau+L}_{\tau}$$, \ie}{}
\[
(I-\Pi_{\cS_n})T(x)=T(\varepsilon),
\]
and by \eqref{tr21},
\begin{equation*}
\begin{aligned}
\left\|(I-\Pi_{\cS_n})T(x)\right\|_\F^2
= \left\| T(\varepsilon)\right\|_\F^2
&= \sum_{|\tau|\leq m} \left\|\Delta^{\tau}\varepsilon\right\|_{n,2}^2
\leq (2m+1)\sigma^2\varkappa^2.
\end{aligned}
\end{equation*}
Due to \eqref{trq1} we conclude that %$\forall\alpha\in(0,1]$,
\[
\left\|T(x)^\H(I-\Pi_{\cS_n})[\zeta]_{-n}^n\right\|_2^2
\le 2(2m+1)\sigma^2\varkappa^2\big(1+\sqrt{\log[1/\alpha]}\big)^2
\]
with probability at least $1-\alpha$. Since
\[\left|\left\langle[\zeta]_{-n}^n, (I-\Pi_{\cS_n})T(x)[\wh\varphi]_{-m}^m\right\rangle\right|\leq \frac{\bar\vrho}{\sqrt{2m+1}} \left\|T(x)^\H(I-\Pi_{\cS_n})[\zeta]_{-n}^n\right\|_{2},
\]
we arrive at the bound with probability $1-\alpha$:
\[
\left|\left\langle[\zeta]_{-n}^n, (I-\Pi_{\cS_n})T(x)[\wh\varphi]_{-m}^m\right\rangle\right|\le\sqrt{2}\sigma\varkappa\bar{\vrho} \big(1+\sqrt{\log[1/\alpha]}\big).
\]
Along with \eqref{i22-1} this results in the bound
\be
%\forall\alpha\in(0,1],\;\;\;\prob\left\{
\prob\left\{\big|{\delta^{(1)}_2}\big|
\leq \sqrt{2}\sigma^2\varkappa(\bar{\vrho}+1) \big(1+\sqrt{\log\left[1/\min (\alpha_2,\alpha_3)\right]}\big)\right\} \ge 1-\alpha_2-\alpha_3.
%\right\}\leq \alpha.
\ee{i22}

\paragraph{3$^o$.} Let us rewrite ${\delta^{(1)}_3}$ as follows:
\bse
{\delta^{(1)}_3} &=&\sigma^2\Re\langle [\zeta]_{-n}^n, (I-\Pi_{\cS_n})M(\wh\varphi){[\zeta]}_{-m-n}^{m+n} \rangle =
\sigma^2 \Re \sigma^2\langle [\zeta]_{-m-n}^{m+n}, Q M(\wh\varphi){[\zeta]}_{-m-n}^{m+n}\rangle,%=\sigma\delta(\zeta,\wh\varphi),
\ese
where $M(\wh\varphi)\in \C^{(2n+1)\times (2m+2n+1)}$ is defined by \eqref{M(phi)}, and $Q\in \C^{(2m+2n+1)\times (2n+1)}$ is given by %$Q_{m,n}=P_{m,n}-E_{m,n}$, where for
\begin{equation*}
%\label{eq:rect_projectors}
Q=
[%\begin{bmatrix}
O_{m,2n+1} ;
I - \Pi_{\cS_n};
O_{m,2n+1}
]%\end{bmatrix}.
%\left[
%\begin{array}{c c}
%0_{m,n+1}; & I - \Pi_{\cS_n}
%\end{array}
%\right].
\end{equation*}
(Hereafter we denote $O_{m,n}$ the $m\times n$ zero matrix.)
Now, by the definition of $\widehat\varphi$ and since the mapping $\varphi \mapsto M(\varphi)$ is linear,
\begin{equation}
\label{maxR}
\begin{aligned}
{\delta^{(1)}_3}
&=\frac{\sigma^2}{2}({[\zeta]}_{-m-n}^{m+n})^\H(\underbrace{Q^{\vphantom{*}}M(\wh\varphi)+M(\wh\varphi)^\H Q^\H}_{K_{1}(\wh\varphi)}) [\zeta]_{-m-n}^{m+n}\\
&\leq\frac{\sigma^2\bar\vrho}{2\sqrt{2m+1}}\max\limits_{\scriptsize\begin{array}{c} u \in \C_m(\Z),\nn \|u\|^*_{m,1}\leq 1\end{array}} ([\zeta]_{-m}^{n})^\H K_{1}(u) [\zeta]_{-m-n}^{m+n}\\
&=\frac{\sigma^2\bar\vrho}{\sqrt{2m+1}}\,\,\max\limits_{\tiny{|j| \leq m}}\,\,\max\limits_{\theta \in [0,2\pi]} \half ([\zeta]_{-m-n}^{m+n})^\H K_{1}(e^{i\theta} u^j) [\zeta]_{-m-n}^{m+n},
%\end{array}
\end{aligned}
\end{equation}
where $u^j \in \C_m(\Z)$, and $[u^j]_{-m}^m = F_m^{\H} e^j$, $e^j$ being the $j$-th canonic basis vector.
Indeed, $([\zeta]_{-m-n}^{m+n})^\H K_{1}(u) [\zeta]_{-m-n}^{m+n}$ is clearly a convex function of the argument $u$ as a linear function of $[\Re(u); \Im(u)]$; as such, it attains its maximum over the set
\begin{align}\label{eq:hyperoct}
\mathcal{B}_{m,1} = \{u \in \C_m(\Z): \|u\|^*_{m,1} \le 1\}
\end{align}
at one of the extremal points~$e^{i\theta} u^j$, $\theta \in [0,2\pi]$, of this set.
It can be directly verified that
\[
K_1(e^{\imath\theta}u) = K_1(u)\cos\theta  + K_2(u)\sin\theta,
\]
where the Hermitian matrix $K_2(u)$ is given by
\begin{align*}
K_2(u) = i \left(Q^{\vphantom{*}}M(u) - M(u)^\H Q^\H \right).
\end{align*}
Denoting $q^j_l(\zeta) = \half([\zeta]_{-m-n}^{m+n})^\H K_l(u^j) [\zeta]_{-m-n}^{m+n}$ for $l=1,2$, we have
\be
\lefteqn{
\max\limits_{\theta \in [0,2\pi]} \half ([\zeta]_{-m-n}^{m+n})^\H K_{1}(e^{\imath\theta} u^j) [\zeta]_{-m-n}^{m+n}
= \max\limits_{\theta \in [0,2\pi]} q^j_1(\zeta)\cos\theta   + q^j_2(\zeta)\sin\theta }\nn
&=&\sqrt{|q_1^j(\zeta)|^2 + |q_2^j(\zeta)|^2} \le \sqrt{2}\max (|q_1^j(\zeta)|, |q_2^j(\zeta)|).
\ee{eq:max_complex_to_real}
Using \eqref{tr22}, by simple algebra we get for $l=1,2$: %for $j=1,...,T+1$,
\[%\begin{align*}
%\label{trA1}
\Tr [K_{l}(u^j)^2]
\leq %\Tr[Q_{m,n}M(u^j)+M(u^j)^*Q_{m,n}^*]
%\Tr([L_T(u[j])+L_T(u[j])^*]^2)
4\Tr[M(u^j)M(u^j)^\H]
= 4(2n+1)\|u^j\|^2_{m,2}
\leq 4(2n+1).
\]%\end{align*}
Now let us bound  $\Tr[K_{l}(u)]$, $l=1,2$, on the set , $\mathcal{B}_{m,1}$ cf.~\eqref{eq:hyperoct}.
One can verify that for the circulant matrix $C(u)$, cf. \eqref{C(phi)}, it holds:
\bse
QM(u) = {R} C(u),
\ese
where $R = Q^{\vphantom{H}}Q^\H$ is an $(2m+2n+1)\times (2m+2n+1)$ projection matrix of rank $s$ defined by
\begin{align*}
R = \left[
\begin{array}{l|l|l}
O_{m,m} & O_{m,n+1} & O_{m,m}\\
\hline
O_{n+1,m} & I-\Pi_{\mathcal{S}_n} & O_{n+1,m} \\\hline
O_{m,m} & O_{m,n+1} & O_{m,m}.
\end{array}
\right]
\end{align*}
Hence, we can bound $\Tr[K_{l}(u)]$, $l=1,2$, as follows:
\begin{align}
|\Tr[K_{l}(u)]|
&\le 2\big|\Tr[R C(u)]\big|\le 2\|R\|_*  \left\|C(u)\right\| \nn
&\le 2\|C(u)\|= 2\sqrt{2m+2n+1}\|\tilde u\|^*_{m+n, 1},
\label{eq:for-zero-padding}
\end{align}
%where $\|u\|^*_{m+n, 1}$ is the Fourier norm of the zero padded
where in the last transition we used the Fourier diagonalization property~\eqref{eq:circ_diag}.
Recall that $u \in \C_m(\Z)$, hence $F_{m+n}[u]$ is the Discrete Fourier Transform of the \textit{zero-padded filter}
\[
\tilde u = [0; ...; 0; [u]_{-m}^m; 0; ... ; 0] \in \C^{2m+2n+1}.
\]
Now combining Lemma \ref{lemma:zero_padding} with~\eqref{eq:for-zero-padding} we arrive at
\[
\left|\Tr[K_{l}(u^j)]\right| \le 2 \sqrt{2m+1}(\kappa_{m,n}^2 + 1)  (\log[2m+2n+1]+3), \quad l=1,2.
\]
By \eqref{trq2} we conclude that for any fixed pair $(l,j) \in \{1,2\} \times \{-m,...,m\}$, with probability $\ge 1-\alpha$,
\[
\big|q^j_l(\zeta)\big| \leq \left|\Tr[K_{l}(u^j)]\right| +\left\|K_{l}(u^j)\right\|_\F \big(1 + \sqrt{\log[2/\alpha]}\big)^2.
\]
With $\alpha_0 = 2(2m+1)\alpha$,  by the union bound together with~\eqref{maxR}~and~\eqref{eq:max_complex_to_real} we get
\be
&\Prob\left\{{\delta^{(1)}_3}
\leq 2\sqrt{2}\sigma^2\bar{\vrho}\left[ (\kappa_{m,n}^2+1)(\log[2m+2n+1]+3)\right.\right.\nn&
\left.\left.\qquad\qquad\qquad\qquad\qquad\qquad\qquad+\kappa_{m,n}\big(1+\sqrt{\log\left[4(2m+1)/ {\alpha_0}\right]}\big)^2\right]\right\}\ge 1-\alpha_0.
\ee{i23-2}

\paragraph{4$^o$.} Bounding ${\delta^{(2)}}$ is relatively easy since $\varphi^{o}$ does not depend on the noise. We decompose
\[
{\delta^{(2)}} = \sigma\Re\langle \zeta, x-\varphi^{o}*x\rangle_n - \sigma^2\Re\langle \zeta,\varphi^{o}*\zeta\rangle_n.
\]
Note that $\Re\langle \zeta, x-\varphi^{o}*x\rangle_n \sim \N(0,\|x-\varphi^{o}* x\|_{n,2}^2)$, therefore, with probability $\ge 1-\alpha$,
\be
\Re\langle\zeta,x-\varphi^{o}*x\rangle_n\leq\sqrt{2\log[1/\alpha]}  \|x-\varphi^{o}*x\|_{n,2}.
\ee{simple0}
On the other hand, defining
\[
\vrho = \sqrt{2m+1} \|\vphi^o\|_{m,1}^*,
\]
we have
\begin{align}\label{simple1}
\|x-\varphi^{o}*x\|_{n,2}
&\leq \|x-\varphi^{o}*y\|_{n,2}+\sigma \|\varphi^{o}*\zeta\|_{n,2}\notag\\
&\leq \|x-\varphi^{o}*y\|_{n,2}+\sqrt{2} \sigma \varrho \kappa_{m,n} \big(1+\sqrt{\log [1/\alpha]}\big)
\end{align}
with probability $1-\alpha$. Indeed, one has
\[
\|\varphi^{o}*\zeta\|^2_{n,2}=\left\|M(\varphi^{o})[\zeta]_{-m-n}^{m+n} \right\|_2^2,
\]
where for $M(\varphi^{o})$ by $\eqref{tr22}$ we have
\begin{equation}
\label{eq:trmm}
\begin{aligned}
\left\|M(\varphi^{o})\right\|_\F^2
= (2n+1)\|\varphi^{o}\|_{m,2}^2
\leq\kappa_{m,n}^2\varrho^2.
\end{aligned}
\end{equation}
Using~\eqref{trq1} we conclude that, with probability at least $1-\alpha$,
\begin{equation}\label{eq:phizeta2}
\|\varphi^{o}*\zeta\|^2_{n,2}
\leq 2\kappa_{m,n}^2\varrho^2\big(1+\sqrt{\log[1/\alpha]}\big)^2,
\end{equation}
which implies~\eqref{simple1}.
Using~\eqref{simple0}~and~\eqref{simple1}, we get that with probability at least $1-\alpha_4-\alpha_5$,
\begin{equation}
\begin{aligned}
&\Re\langle\zeta,x-\varphi^{o}*x\rangle_n \\
&\;\leq \sqrt{2\log \left[1/\min(\alpha_4,\alpha_5)\right]}
\left[\|x-\varphi^{o}*y\|_{n,2}+\sqrt{2}\sigma\varrho\kappa_{m,n} \big(1+\sqrt{\log[1/\min(\alpha_4,\alpha_5)]}\big)\right]\\
&\;\leq \|x-\varphi^{o}*y\|_{n,2}\sqrt{2\log\left[1/\min(\alpha_4,\alpha_5)\right]}+2\sigma\varrho \kappa_{m,n}\big(1+\sqrt{\log \left[1/\min(\alpha_4,\alpha_5)\right]}\big)^2.
\end{aligned}
\label{simple2}
\end{equation}
Now, the (indefinite) quadratic form
\begin{align*}
&\Re\langle\zeta,\varphi^{o}*\zeta\rangle_n=\half ([\zeta]_{-m-n}^{m+n})^\H{K_0(\varphi^{o})}[\zeta]_{-m-n}^{m+n},
\end{align*}
where
\[
K_0(\varphi^{o}) =
[
O_{m,2m+2n+1};
M(\varphi^{o});
O_{m,2m+2n+1}]
+
[O_{m,2m+2n+1} ;
M(\varphi^{o});
O_{m,2m+2n+1}
]^\H,
\]
%can be bounded as in {\bf Step 3$^o$}, so that
whence (cf. {\bf 3$^o$})
\begin{align*}
|\Tr [{K_0(\varphi^{o})}]|
&\le 2(2n+1)\left|\varphi^{o}_0\right|
%&\leq 2\kappa_{m,n}^2\varrho.
\end{align*}
Let us bound $\left|\varphi^{o}_0\right|$.
Let~$e^{0}$ be the discrete centered Dirac vector in $\R^{2m+1}$, and note that $\|F_m [{e}^0]\|_{\infty} = 1/\sqrt{2m+1}$.
Then,
\begin{align*}
\left|\varphi^{o}_m\right|
= |\langle [\varphi^{o}]_{-m}^m, {e}^0 \rangle|
\le \|\varphi^{o}\|^*_{m,1} \|F_m [{e}^0]\|_{\infty}
\le \frac{\vrho}{2m+1},
\end{align*}
whence $|\Tr[{K_0(\varphi^{o})}]| \leq 2\kappa_{m,n}^2\vrho$.
%since obviously $\|F_m [{e}^0]\|_{\infty} = 1/\sqrt{2m+1}$.
On the other hand, by \eqref{eq:trmm},
\[
\left\|K_0(\varphi^{o})\right\|_\F^2
\leq 4\left\|M(\varphi^{o})\right\|_\F^2
\le 4\kappa_{m,n}^2\varrho^2.
\]
Hence by \eqref{trq2},
\begin{equation}
\label{i112}
\Prob\left\{-\Re \langle\zeta,\varphi^{o}*\zeta\rangle_n \leq 2\kappa_{m,n}^2\varrho+2\kappa_{m,n}\varrho \big(1+\sqrt{2\log\left[1/\alpha_6\right]}\big)^2 \right\} \ge 1-\alpha_6.
\end{equation}

\paragraph{5$^o$.}
Let us combine the bounds obtained in the previous steps with initial bound~$\eqref{first0}$. For any $\alpha \in (0,1]$, putting $\alpha_i = \alpha/4$ for $i = 0, 1, 6$, and $\alpha_j = \alpha/16$, $2 \le j \le 5$, by the union bound we get that with probability $\ge 1-\alpha$,
\begin{eqnarray}
\|x-\wh\varphi*y\|_{n,2}^2
&\leq& \|x-\varphi^{o}*y\|_{n,2}^2+2{\delta^{(2)}}-2{\delta^{(1)}}\nn
\mbox{[by \eqref{simple2}]}&\leq& \|x-\varphi^{o}*y\|_{n,2}^2+2\sigma\|x-\varphi^{o}*y\|_{n,2}\sqrt{2\log[16/\alpha]}\nn
\mbox{[by \eqref{simple2}--\eqref{i112}]}
&&+ \, 4\sigma^2\varrho\Big[\kappa_{m,n}^2+2\kappa_{m,n}\big(1+\sqrt{2\log[16/\alpha]}\big)^2\Big]
%\label{eq:final_long-1}
\nn
\mbox{[by \eqref{i11}]}
&&+ \, 2\sigma\|x-\wh\varphi*y\|_{n,2}\big(\sqrt{2s}+\sqrt{2\log [16/\alpha]}\big)
%\label{eq:final_long-2}
\nn
\mbox{[by \eqref{i22}]}
&&+ \, 2\sqrt{2}\sigma^2(\bar{\vrho}+1)\big(1+\sqrt{\log[16/\alpha]}\big)\varkappa
%\label{eq:final_long-3}
\nn
\mbox{[by \eqref{i23-2}]}
&&+ \, 4\sqrt{2}\sigma^2\bar{\vrho} \Big[(\kappa_{m,n}^2+1)(\log[2m+2n+1]+3) \nn &&\qquad\qquad\quad+\kappa_{m,n}\left(1+\sqrt{\log\left[{16(m+1)}/{\alpha}\right]}\right)^2\Big]  %\label{eq:final_long-4}
\label{eq:l2con-long}
\end{eqnarray}
Now, denote $c_\alpha := \sqrt{2\log[16/\alpha]}$ and let
\begin{align}
u(\alpha) &= 2\big( \sqrt{2}+c_\alpha\big),\label{eq:part-u}\\
v_1(\alpha) &= 4\left[ \kappa_{m,n}^2+2\kappa_{m,n}\left(1+c_\alpha\right)^2\right],\label{eq:part-v1}\\
v_2(\alpha )&= 4\sqrt{2} \Big[ (\kappa_{m,n}^2+1)(\log[2m+2n+1]+3) +\kappa_{m,n}\left(1+\sqrt{\log\left[{16(2m+1)}/{\alpha}\right]}\right)^2\Big].
\label{eq:part-v2}
\end{align}
In this notation, \eqref{eq:l2con-long} becomes
\begin{align}
\|x-\wh\varphi*y\|_{n,2}^2
\le \|x-\varphi^{o}*y\|_{n,2}^2
&+ 2\sigma (\sqrt{2s}+c_\alpha) \left(\|x-\wh\varphi*y\|_{n,2} + \|x-\varphi^{o}*y\|_{n,2}\right)\nn
&+ u(\alpha)\sigma^2(\bar\varrho+1)\varkappa + (v_1(\alpha)+v_2(\alpha)) \sigma^2\bar\varrho,
\label{eq:uv-compact}
\end{align}
which implies, by completing the squares, that
\begin{align*}
\|x-\wh\varphi*y\|_{n,2} \le \|x-\varphi^{o}*y\|_{n,2}
+2\sigma(\sqrt{2s}+c_\alpha) + \sigma\sqrt{u(\alpha)(\bar\varrho+1)\varkappa + (v_1(\alpha)+v_2(\alpha))\bar\varrho}.
\end{align*}
Let us simplify this bound. Note that
\begin{equation}\label{eq:part-u-bound}
u(\alpha) \le 4c_\alpha,
\end{equation}
while on the other hand,
\begin{align}
\label{eq:final_stoch_term}
v_1(\alpha) + v_2(\alpha)
&\le 4\sqrt{2} (\kappa_{m,n}^2+1)(\log[2m+2n+1]+4) + 4.5(4\sqrt{2}+8)\kappa_{m,n}\log\left[16(2m+1)/\alpha\right]\nn
&\le 8\left(1 + 4\kappa_{m,n}\right)^2 \log\left[110(m+n+1)/\alpha\right].% = 4V_\alpha.
\end{align}
We finally arrive at
\begin{equation}\label{eq:l2con-mn-stat}
\|x - \widehat{\varphi}*y \|_{n,2} \le \|x-\varphi^{o}*y\|_{n,2} + 2 \sigma \left( \sqrt{\bar\varrho V_\alpha} + \sqrt{(\bar{\varrho}+1)c_\alpha \varkappa} + \sqrt{2s}+ c_\alpha \right)
\end{equation}
where we put
\begin{equation}
\label{eq:V}
V_\alpha := 2\left(1 + 4\kappa_{m,n}\right)^2\log \left[{110(m+n+1)}/{\alpha}\right].
\end{equation}
The bound~\eqref{eq:th1} of the theorem follows from~\eqref{eq:l2con-mn-stat} after straightforward simplifications.
\qed

\subsection{Proof of Theorem \ref{th:l2pen++}}
%Let $\widehat\varrho = \sqrt{n+1} \cdot \|\widehat\varphi\|^*_{n,1}$.
Denote~$\wh\vrho = \sqrt{2m+1}\|\wh\vphi\|_{m,1}^*$, and let~$\vrho = \vrho(\vphi^o)=\sqrt{2m+1}\|\vphi^o\|_{m,1}^*$ for some $\vphi^o\in \C_m(\Z)$.
In the sequel, we use the notation defined in the proof of Theorem \ref{th:l2con}.
%with $\widehat\varrho := \sqrt{m+1}\,\|\widehat\varphi\|^*_{m,1}$,
We have the following counterpart of~\eqref{first0}:
\begin{equation*}
\|x-\wh\varphi*y\|_{n,2}^2 + \lambda^2\sigma^2\widehat\varrho^2
\leq \|x-\varphi^{o}*y\|_{n,2}^2 - 2\delta^{(1)} + 2\delta^{(2)} + \lambda^2\sigma^2\varrho^2.
\end{equation*}
When repeating steps {\bf 1$^o$--4$^o$} of the proof of Theorem \ref{th:l2con} we obtain a counterpart of~\eqref{eq:uv-compact}:
\be
\|x-\wh\varphi*y\|_{n,2}^2 + \lambda^2\sigma^2\widehat\varrho^2
&\leq
&\|x-\varphi^{o}*y\|_{n,2}^2 + 2\sigma (\|x-\varphi^{o}*y\|_{n,2} \nn&&+ \|x-\wh\varphi*y\|_{n,2})(\sqrt{2s}+c_\alpha) + u(\alpha)\sigma^2\varkappa+v_1(\alpha)\sigma^2\varrho\nn
&& + \lambda^2\sigma^2\varrho^2 + \left[u(\alpha)\varkappa+v_2(\alpha)\right]\sigma^2\wh\varrho
\ee{eq:l23compete}
with~$u(\alpha)$, $v_1(\alpha)$, and~$v_2(\alpha)$ given by~\eqref{eq:part-u}--\eqref{eq:part-v2}.
We now consider two cases as follows.
\paragraph{(a)}
First, assume that
\be
\|x-\wh\varphi*y\|_{n,2}^2 &\leq& \|x-\varphi^{o}*y\|_{n,2}^2 + 2\sigma (\|x-\varphi^{o}*y\|_{n,2} + \|x-\wh\varphi*y\|_{n,2}) (\sqrt{2s}+c_\alpha) \nn
&&\qquad\qquad\qquad\quad+ u(\alpha)\sigma^2\varkappa+v_1(\alpha)\sigma^2\varrho + \lambda^2\sigma^2\varrho^2.
\ee{eq:assume1}
In this case, clearly,
\be
\|x-\wh\varphi*y\|_{n,2} &\leq & \|x-\varphi^{o}*y\|_{n,2}+2\sigma\big(\sqrt{2s}+c_\alpha\big)
+\sqrt{u(\alpha)\sigma^2\varkappa+v_1(\alpha)\sigma^2\varrho
 + \lambda^2\sigma^2\varrho^2}\nn
 &\leq &\|x-\varphi^{o}*y\|_{n,2} + 2\sigma (\sqrt{2s}+c_\alpha)+\sigma(\sqrt{u(\alpha)\varkappa+v_1(\alpha)\varrho}
 + \lambda\varrho)
\ee{eq:1stcase}

\paragraph{(b)}
Suppose, on the contrary, that~\eqref{eq:assume1} does not hold, we then conclude from~\eqref{eq:l23compete} that
%\[
%\lambda^2\wh{\varrho}^2 \leq \wh\varrho(u(\alpha)\varkappa+v_2(\alpha)),
%\]
%same as
\[
%\wh\varrho \leq 2\sqrt{2}\sigma^2\lambda^{-2}\overbrace{\left[\varkappa\big(1+\sqrt{\log[16/\alpha]}\big)+{2\left[(\kappa_{m,n}^2+1)(\log[m+n+1]+3)
%+\kappa_{m,n}\left(1+\sqrt{\log\left[{16(m+1)}/{\alpha}\right]}\right)^2\right]}\right]}^{\varkappa\big(1+\sqrt{\log[16/\alpha]}\big)
%+\varsigma_1(\alpha)}
\wh\varrho \leq\lambda^{-2}(u(\alpha)\varkappa+v_2(\alpha)),
\]
and
\[
u(\alpha)\wh\varrho\varkappa+v_2(\alpha)\wh\varrho
\leq\lambda^{-2}(u(\alpha)\varkappa+v_2(\alpha))^2.
\]
When substituting the latter bound into~\eqref{eq:l23compete}, we obtain the bound
\bse
{\|x-\wh\varphi*y\|_{n,2}}
&\leq& \|x-\varphi^{o}*y\|_{n,2} +2\sigma(\sqrt{2s}+c_\alpha)\\&&
+ \sigma(\sqrt{u(\alpha)\varkappa+v_1(\alpha)\varrho} + \lambda^{-1} (u(\alpha)\varkappa+v_2(\alpha))+ \lambda\varrho),
\ese
which also holds in the case of {\bf (a)} due to~\eqref{eq:1stcase}.\\

Finally, using~\eqref{eq:part-u-bound},~\eqref{eq:final_stoch_term}, and the bound
\bse
v_1(\alpha) \le 4(1+\kappa_{m,n})^2(1+c_\alpha)^2
\ese
which directly follows from~\eqref{eq:part-v1}, we conclude that
\bse
\|x - \widehat{\varphi}*y \|_{n,2}
&\le& \|x-\varphi^{o}*y\|_{n,2} \\&&+\sigma(\lambda \varrho + 4\lambda^{-1}(c_\alpha \varkappa + V_\alpha)) + 2 \sigma \left(\sqrt{\varrho W_\alpha} + \sqrt{c_\alpha \varkappa} + \sqrt{2s} +  c_\alpha \right)
\ese
with $V_\alpha$  given by~\eqref{eq:V}, and
$
W_\alpha = (1+\kappa_{m,n})^2 (1+c_\alpha)^2.
$
The bound~\eqref{eq:l2pen++} of the theorem follows by a straightforward simplification of the above bound.
\qed

%\section{Proofs for Section \ref{sec:pointwise}}\label{sec:prpointwise}

\section{Proofs for Section \ref{sec:shift-inv}}

\subsection{Proof of Proposition~\ref{prop:shift-inv-bilateral}}
%Assume that $x$ holds for some $m$ such that  $m+1\geq s$ and $\varkappa = 0$.
Let $\Pi_{\S_m}$ be the $m+1$-dimensional Euclidean projection matrix on the subspace $\S_m\subset \C^{m+1}$ of dimension $\leq s$ (in fact, this subspace is exactly of dimension $s$) generated by vectors $x_{0}^m$ for $x\in \S$ (one may set, for instance,  $\Pi_{\S_m}=Z_m(Z_m^HZ_m)^{-1}Z^H_m$, $Z_m=[z_1,...,z_{\dim(\S_m)}]$, where $z_i$ are linearly independent and such that $z_i=[x_i]_0^m$ with $x_i\in \S$).
Since $\dim(\S)\leq s$, one has
\[
\|\Pi_{\S_m}\|_2^2 = \Tr(\Pi_{\S_m}) \leq s.
\]
Thus, there is a $j \in \{0,...,m\}$ such that the $j+1$-th {column} $r = [\Pi_{\S_m}]_j$ of~$\Pi_{\S_m}$ satisfies
\[
\left\|r\right\|_2\leq \sqrt{\frac{s}{m+1}} \le \sqrt{\frac{2s}{2m+1}},
\]
and, because $\Pi_{\S_m}$ is the projector on $\S_m$ one has $x_{j}-\langle r , x_0^m\rangle=0$ for all $x \in \S$.
Hence, using that $\Delta \S = \S$ we obtain for all $\tau \in \Z$
\[
x_\tau-\langle r, x_{\tau-j}^{\tau-j+m}\rangle=0, \quad \tau \in \Z.
\]
%It remains to augment $\pi$ with zeroes in such a way that the augmented vector corresponds to $\phi^o \in \C_m(\Z)$.
Finally, let $\phi^o \in \C_m(\Z)=\Delta^{-j}\phi(r)$ where $\phi(r)$ is the inverse slicing map of $\tilde r\in \C^{m+1}$ such that $\tilde r_i=r_{m+1-i}$. Obviously, $\phi^o \in \C_m(\Z)$; on the other hand,
\[
\|\phi^o\|_2\leq \sqrt{\frac{2s}{2m+1}} \quad \text{and} \quad x_t-[\phi^o*x]_t=0, \quad \forall t \in \Z. \tag*{\qed}
\end{equation*}

\subsection{Proof of Proposition \ref{pr:l2sloss}}
\label{sec:prmisc}
%When proving the results of Section \ref{sec:shift-inv}
In the proofs to follow, the following simple statement will be of use.
\begin{lemma}
\label{le:shift-inv-approx-err}
\item[(i)] Suppose that for all $z \in \S$ there is a filter $\phi^o\in \C_m(\Z)$ such that $z=\phi^o* z$ with $\|\phi^o\|_2\leq{\rho\over \sqrt{2m+1}}$ for some $\rho \ge 1$.
Then for all $x\in\X_{m,n}(s,\varkappa)$ one has
\be
\|x-\phi^o  *x\|_{n,2}\leq \sigma\varkappa(1+\rho).
\ee{eq:s0001}
Moreover, if $x\in \overline\X_{m,n}(s,\varkappa)$ then
\be
\|x-\phi^o*x\|_{n,\infty}\leq {\sigma\varkappa \over \sqrt{2m+1}}(1+\rho\kappa_{n,m}).
\ee{eq:s0002}
\item[(ii)] Similarly, assume that for all $z \in \S$ there is $\phi^o\in \C_m(\Z)$ such that $z=\phi^o* \Delta^m z$ and $\|\phi^o\|_2\leq {\rho\over \sqrt{2m+1}}$ for some $\rho \ge 1$.
Then for all $x\in\X_{m,n}(s,\varkappa)$ one has
\[
\|\Delta^{-m}(x-\phi^o*\Delta^mx)\|_{n,2}\leq \sigma\varkappa(1+\rho).
\]
Furthermore, if $x\in \overline\X_{m,n}(s,\varkappa)$ then
\[
\|\Delta^{-m}(x-\phi^o*\Delta^mx)\|_{n,\infty}\leq {\sigma\varkappa \over \sqrt{2m+1}}(1+\rho\kappa_{n,m}).
\]
\end{lemma}
\paragraph{Proof of the lemma.} Here we prove the first statement of the lemma, proof of the second one being completely analogous. Recall that any $x  \in\X_{m,n}(s,\varkappa)$ can be decomposed as in
$x=x^\S+\varepsilon$ where $x^\S\in \S$ and $\|\Delta^\tau\varepsilon\|_{n,2}\leq \varkappa \sigma$ for all $|\tau|\leq m$. Thus,
\be
\|x - [\phi^o * x]\|_{n,2}\le \|x^\S - \phi^o * x^\S\|_{n,2} + \|\varepsilon\|_{n,2} + \|\phi^o * \varepsilon\|_{n,2}=
\varkappa \sigma+\|\phi^o * \varepsilon\|_{n,2}.
\ee{eq:p0001}
On the other hand, by the Cauchy inequality,
\[\|\phi^o * \varepsilon\|^2_{n,2}=\sum_{t=-n}^n\left|
\sum_{\tau=-m}^m\phi^o_\tau\varepsilon_{t-\tau}\right|^2
\leq \|\phi^o\|^2_{2}\sum_{t=-n}^n\sum_{\tau=-m}^m|\varepsilon_{t-\tau}|^2=\|\phi^o\|^2_{2}\sum_{\tau=-m}^m\|\Delta^\tau\varepsilon\|^2_{n,2}\leq
\rho^2\sigma^2\varkappa^2.
\]
When substituting the latter bound into \rf{eq:p0001} we obtain \rf{eq:s0001}.
\par
To show  \rf{eq:s0002} recall that in the case of $x\in \overline \X_{m,n}(s,\varkappa)$ we have
$x=x^\S+\varepsilon$ with $|\varepsilon_{\tau}|\leq {\varkappa \sigma\over \sqrt{2n+1}}$ for all $|\tau|\leq m+n$. Then for $|t|\leq n$ we get \bse
|x_t - [\phi^o * x]_t|&\le& |x^\S_t - [\phi^o * x^\S]_t| + |\varepsilon_t| + |[\phi^o * \varepsilon]_t|\leq
 {\varkappa \sigma\over \sqrt{2n+1}}+\|\phi^o\|_2\|\Delta^{-t} \varepsilon\|_{m,2}\\&\leq& {\varkappa \sigma\over \sqrt{2n+1}} +{\rho\over \sqrt{2m+1}} {\sigma\varkappa\sqrt{2m+1} \over \sqrt{2n+1}}\leq {\sigma\varkappa \over \sqrt{2m+1}}(1+\rho\kappa_{n,m}).\qquad\qquad\mbox{\qed}
\ese
\paragraph{Proof of the proposition.} W.l.o.g.~we may assume that $m=2m_o$.
In the premise of the proposition, by Proposition \ref{prop:shift-inv-bilateral}, for any ${m_o}\geq s-1$ 
{there exists} a filter $\phi^o\in \C_{m_o}(\Z)$ such that
\be
\|\phi_o\|_2\leq \sqrt{2s\over 2m_o+1},\;\;z=\phi_o*z\;\;\forall z\in \S.
\ee{eq:{m_o}f}
When setting $\vphi^o=\phi^o*\phi^o\in \C_m$ we have $z-\vphi^o*z=0$ $\forall z\in \S$, and\footnote{In the case of $m=2m_o+1$ one may consider two filters $\phi^o$ and $\psi^o$ of widths $m_o$ and $m_o+1$ respectively, and then build  $\vphi_o=\phi^o *\psi^o\in \C_m(\Z)$. One easily verifies that in this case
$\|\vphi_o\|^*_{m,1}\leq \sqrt{2m+1}\|\phi^o\|_2\|\psi^o\|_2\leq {4s\over \sqrt{2m+1}}$.}
\be
\|\vphi^o\|_{m,2}\leq \|\vphi^o\|^*_{m,1}&\leq& {4s\over \sqrt{2m+1}}
\ee{eq:nof*2}
(cf. \cite[Proposition 3]{harchaoui2015adaptive} or \cite[Lemma 16]{jn1-2009}).
We now apply Lemma \ref{le:shift-inv-approx-err}.i to obtain for all $x\in  \X_{m,n}(s,\varkappa)$
\be
\|x-\vphi^o*x\|_{n,2}\leq \sigma\varkappa (4s+1).
\ee{eq:obias}
%We need the following result (cf. \cite[Proposition 3]{harchaoui2015adaptive} and \cite[Lemma 16]{jn1-2009}).
%Let $m_o\in \Z_+$, $m=2m_o$, and $\phi^o\in \C_{m_o}(\Z)$ such that $\|\phi^o\|_2\leq {\rho\over \sqrt{2m_o+1}}$ and
%\[
%\|x-\phi^o*x\|_{n+m_o,2}\leq \sigma\theta.
%\]
%Then the filter $\vphi^o=\phi^o*\phi^o\in \C_m$ satisfies
%\be
%\|\vphi^o\|^*_{m,1}&\leq& {2\rho^2\over \sqrt{2m+1}},\;\;\;\|\vphi^o\|_{1}\leq \rho^2,\nn%\;\mbox{and}\;
%\|x-\vphi^o*x\|_{n,2}&\leq& \sqrt{2}\kappa_{m,n}\sigma\theta(\rho+1).
%\ee{eq:obias}
Moreover, note that
\[
\|\varphi^o*\zeta\|_{n,2}^2=\langle \zeta,M(\varphi^o)\zeta\rangle_{n},
\]
where $M(\varphi)$ is defined by~\eqref{M(phi)}. When using the bound \rf{eq:nof*2} along with \eqref{tr22} we obtain
\[
\|M(\varphi^o)\|_\F^2=(2n+1)\|\varphi^o\|_2^2\leq 16\kappa_{m,n}^2s^2;
\]
by \eqref{trq1} this implies that for any $\alpha\in(0,1)$, with probability at least $1-\alpha$,
\be
\|\varphi^o*\zeta\|_{n,2}\leq 4\sqrt{2}\sigma\kappa_{m,n} s\big(1+\sqrt{\log[1/\alpha]}\big).
\ee{eq:varphizeta}
The latter bound taken together with \rf{eq:obias}  implies that with probability $\geq 1-\alpha$
\bse
\|x-\vphi^o*y\|_{n,2}&\leq& 4\sqrt{2}\kappa_{m,n}\sigma s\big(1+\sqrt{\log[1/\alpha]}\big)+\sigma\varkappa(4s+1)\\
&\leq &
C\sigma s\big(\kappa_{m,n}\sqrt{\log[1/\alpha]}+\varkappa\big)
\ese
when $\alpha\leq 1/2$. We conclude the proof by substituting the above bound for the loss of the estimate $\wh x=\vphi^o*y$ and the bound $\|\vphi^o\|^*_{m,1}\leq 4s$ into the oracle inequalities of Theorems \ref{th:l2con} and \ref{th:l2pen++}.
\qed

\subsection{Proof of Proposition \ref{pr:linfsloss}}
We provide the proof for the case of constrained estimator $\wh x_{\con}$, the proof of the proposition for penalized estimator $\wh x_{\pen}$ follows exactly same lines. Let $\wh\vphi=\wh\vphi_{\con}$.
\paragraph{1$^o$.}W.l.o.g. we assume that $m=2m_o$.
By Proposition \ref{prop:shift-inv-bilateral}, for such ${m_o}$ there is a filter $\phi^o\in \C_{m_o}(\Z)$ satisfying relationships \rf{eq:{m_o}f}.
When applying Lemma \ref{le:shift-inv-approx-err}.i we obtain for all $x\in{\overline \X_{m,n}(s,\varkappa)}$
\be
\|x-\phi^o*x\|_{n,\infty}\leq {\sigma\varkappa \over \sqrt{2m_o+1}}(1+\sqrt{2s}\kappa_{n,m_o}).
%\leq {2\sigma\varkappa \over \sqrt{2m+1}}(1+\sqrt{2s}\kappa_{m,n}).
\ee{eq:obiasinf}
Next, replacing $\vphi^o$ with $\phi^o$ and $n$ with $m$ in the derivation which led us to \rf{eq:varphizeta} in the proof of Proposition \ref{pr:l2sloss} we conclude that
\be
\|\phi^o*\zeta\|_{m,2}\leq 2\sigma\kappa_{m_o,m} \sqrt{s}\big(1+\sqrt{\log[1/\alpha]}\big)\leq 2\sqrt{2s}\sigma \big(1+\sqrt{\log[1/\alpha]}\big).
\ee{eq:ovarinf}
\paragraph{2$^o$.} Let now $|t|\leq n-m_o$. We decompose
\be
|[x - \widehat\varphi * y]_t|
&=&|[(\phi^o+(1-\phi^o))*(x - \widehat\varphi * y)]_t|\nn
&\leq& |[\phi^o*(x - \widehat\varphi * y)]_t| + |[(1 - \phi^o) *(1-\widehat \varphi) *  x]_t|\nn
&&+ \sigma|[\widehat \varphi * \zeta]_t| + \sigma|[\widehat \varphi * \phi^o * \zeta]_t|\nn
&=:&\delta^{(1)}+\delta^{(2)}+\delta^{(3)}+\delta^{(4)}.
\ee{eq:l2point-long}
We have
\[
\delta^{(1)}
\leq \|\phi^o\|_2 \big\|\Delta^{-t} [x - \widehat\varphi * y]\big\|_{{m_o},2} \nn
\leq {2\sqrt{s}\over \sqrt{2m+1}} \|x - \widehat\varphi * y\|_{n,2}.
\]
%where we used that $n \le t \le 2n$.
Using the bound \rf{eq:riskl2-2} of Proposition \ref{pr:l2sloss}  we conclude that with probability $\ge 1-\alpha/3$
\[
\delta^{(1)} \leq  C{\sqrt{s\over 2m+1}}\overline\psi^{\alpha/3}_{m,n}(\sigma,s,\varkappa).
\]
%\paragraph{2$^o$.}
Next, using \rf{eq:obiasinf} we get
\[
\delta^{(2)}
\le \left(1 + \left\|\widehat\varphi\right\|_1\right) \left\|\Delta^{-t}[(1 - \phi^o) * x]\right\|_{{m_o},\infty}
\le C's{\sigma\varkappa \over \sqrt{2m+1}}(1+\sqrt{2s}\kappa_{n,m_o})\le {Cs^{3/2}\sigma\varkappa\over\sqrt{2m+1}}
\]
(recall that $n\geq m_o$).
%where the last transition is due to $n \ge {m_o}$.
%where we again used that $n \le t \le 2n$.
%\paragraph{3$^o$.}
Further, by the Parseval's identity, with probability $\geq 1-\alpha/3$,
\[
\delta^{(3)}
=\sigma|\langle F_{m}[\wh\vphi], F_{m}[\Delta^{-t} \zeta]\rangle|
\leq \sigma\|\widehat{\varphi}\|^*_{m,1} \|\Delta^{-t}\zeta\|^*_{m,\infty}\le \frac{C's\sigma}{\sqrt{2m+1}}\sqrt{2 \log \left[3(2m+1)/\alpha\right]}
\]
 due to \eqref{eq:complex_gaussian_max}.
%\par
Finally, using \rf{eq:ovarinf} and the fact that the distribution of $\zeta_{t-m-m_o}^{t+m+m_o}$ is the same as that of $\zeta_{-m-m_o}^{m+m_o}$ we conclude that with probability $\geq 1-\alpha/3$ it holds
\[
\|\Delta^{-t}[\phi^o * \zeta]\|_{m,2}  \leq 2\sqrt{2s}\sigma \big(1+\sqrt{\log[3/\alpha]}\big).
\]
Therefore, we have for $\delta^{(4)}$:
\bse
\delta^{(4)}
&\le& \sigma\|\widehat \varphi\|_{m,2} \left\| \Delta^{-t} [\phi^o * \zeta] \right\|_{m,2}
\le \frac{C's\sigma}{\sqrt{2m+1}} 2\sqrt{2s}\sigma \big(1+\sqrt{\log[3/\alpha]}\big)\nn
&=& \frac{C'' s^{3/2} \sigma}{\sqrt{2m+1}}\left(1+\sqrt{\log[3/\alpha]}\right)
\ese
with prob.~$\ge 1-\alpha/3$.
Substituting the bounds for $\delta^{(k)},\,k=1,...,4$, into \eqref{eq:l2point-long} we arrive at \rf{eq:pi}. \qed

\subsection{Proof of Proposition \ref{th:shift_invariant}}

%As discussed in Remark \ref{rem:diff_eq}, there is a one-to-one correspondence between shift invariant subspaces of $\C(\Z)$ and homogeneous difference equations with polynomial operators.
As a precursory remark, note that if a finite-dimensional subspace $\S$ is shift-invariant, i.e.,~$\Delta \S \subseteq \S$, then necessarily $\Delta\S = \S$ (indeed, $\Delta$ obviously is a linear transformation with a trivial kernel).

\paragraph{1$^o$.}
To prove the direct statement, note that the solution set of \eqref{eq:shift-inv-diff_eq} with $\deg({p}(\cdot)) = s$ is a shift-invariant subspace of $\C(\Z)$ -- let us call it $\S'$. Indeed, if $x \in \C(\Z)$ satisfies~\eqref{eq:shift-inv-diff_eq}, so does $\Delta x$, so $\S'$ is shift-invariant. To see that $\dim(\S') = s$, note that $x \mapsto x_1^s$ is a bijection $\S' \to \C^s$: under this map arbitrary $x_1^s \in \C^{s}$ has a unique preimage. Indeed, as soon as one fixes $x_1^s$, \eqref{eq:shift-inv-diff_eq} uniquely defines the next samples $x_{s+1}, x_{s+2}, ...$ (note that ${p}(0) \ne 0$); dividing \eqref{eq:shift-inv-diff_eq} by $\Delta^s$, one can retrieve the remaining samples of $x$ since $\deg({p}(\cdot))=s$ (we used that $\Delta$ is bijective on $\S$).

\paragraph{2$^o$.}
To prove the {converse,} first note that any polynomial ${p}(\cdot)$ with $\deg({p}(\cdot)) = s$ and such that ${p}(0)=1$ is uniquely expressed via its roots $z_1, ..., z_s$ as
\[
{p}(z) = \prod_{k=1}^s (1-z/z_k).
\]
Since $\S$ is shift-invariant, we have $\Delta \S = \S$ as discussed above, i.e.,~$\Delta$ is a bijective linear operator on~$\S$. Let us fix some basis $E = [e^1; ...;  e^s]$ of $\S$ and denote $A$ the $s \times s$ representation matrix of $\Delta$ in this basis, that is, $\Delta(e^j) = \sum_{i=1}^s a_{ij} e^{i}$. By the Jordan theorem basis $E$ can be chosen in such a way that $A$ is upper-triangular.
Then, any vector $x \in \S$ satisfies
$%\label{eq:diff_diag}
{q}(\Delta) x \equiv 0$
where
\[
\begin{aligned}
{q}(z)
&= \prod_{i=1}^s (a_{ii}-z) = \det(A-z I)
\end{aligned}
\]
is the characteristic polynomial of $A$.
Note that $\det A = \prod_{i=1}^s a_{ii} \ne 0$ since $\Delta$ is a bijection.
Hence, choosing
\[
{p}(\Delta) = \frac{{q}(\Delta)}{\det A}
\]
we obtain
$\prod_{i=1}^s (1- c_i\Delta) x \equiv 0$ for some complex $c_i \ne 0$.
This means that $\S$ is contained in the solution set $\S'$ of \eqref{eq:shift-inv-diff_eq} with $\deg({p}(\cdot))=s$ and such that ${p}(0)=1$.
Note that by~$\bf{1}^o$ $\S'$ is also a shift-invariant subspace of dimension $s$, thus $\S$ and $\S'$ coincide.
Finally, uniqueness of ${p}(\cdot)$ follows from the fact that ${q}(\cdot)$ is a characteristic polynomial of $A$.
\qed

\subsection{Proof of Proposition~\ref{th:sines}}
To prove the proposition we need to exhibit a vector $q\in \C^{n+1}$ of small $\ell_2$-norm and such that the polynomial
$1-q(z)=1-\left[\sum_{i=0}^n q_iz^i\right]$ is divisible by $p(z)$, i.e.,~that there is a polynomial $r(z)$ of degree $n-s$ such that
\[
1-q(z)=r(z)p(z).
\]
Indeed, this would imply that
\[
x_t-[q*x]_t=[1-q(\Delta)]x_t=r(\Delta)p(\Delta)x_t=0
\]
due to $p(\Delta)x_t=0$,
\par
%The bound
%$
%\|q\|_2\leq C's^{3/2}\sqrt{\log s\over m}
%$
%of \eqref{propofq} is proved in \cite[Lemma 6.1]{jn-2014}.
Our objective is to prove the inequality
$ \|q\|_2\leq C's\sqrt{\log[ns]\over n}.$
So, let $\theta_1,...,\theta_s$ be complex numbers of modulus 1 -- the roots of the polynomial $p(z)$. Given $\delta=1-\epsilon\in(0,1)$, let us set $\bar{\delta}={2\delta/(1+\delta)}$, so that
\begin{equation}\label{eq1}
{\bar{\delta}\over\delta}-1=1-\bar{\delta}>0.
\end{equation}
Consider the function
\[%begin{equation}\label{eq2}
\bar{q}(z)=\prod_{i=1}^s{z-\theta_i\over\delta z-\theta_i}.
\]%end{equation}
\def\B{{\cal B}}
Note that $\bar{q}(\cdot)$ has no singularities in the circle
\[
\B=\{z:|z|\leq {1/\bar{\delta}}\};
\]
besides this, we have
$%\begin{equation}\label{eq4}
\bar{q}(0)=1.
$
Let $|z|=1/\bar{\delta}$, so that $z=\bar{\delta}^{-1}w$ with $|w|=1$. We have
$$
{|z-\theta_i|\over|\delta z-\theta_i|}={1\over\delta}{|w-\bar{\delta}\theta_i|\over |w-{\bar{\delta}\over \delta}\theta_i|}.
$$
 We claim that when $|w|=1$, $|w-\bar{\delta}\theta_i|\leq |w-{\bar{\delta}\over\delta}\theta_i|$.
\begin{quote}
Indeed, assuming w.l.o.g. that $w$ is not proportional to $\theta_i$, consider triangle $\Delta$ with the vertices $A=w$, $B=\bar{\delta}\theta_i$ and $C={\bar{\delta}\over\delta}\theta_i$. Let also $D=\theta_i$. By (\ref{eq1}), the segment $\overline{AD}$ is a median in $\Delta$, and $\angle CDA $ is $\geq{\pi\over 2}$ (since $D$ is the closest to $C$ point in the unit circle, and the latter contains $A$), so that  $|w-\bar{\delta}\theta_i|\leq |w-{\bar{\delta}\over\delta}\theta_i|$.
\end{quote}
As a consequence, we get
\begin{equation}\label{eq10}
z\in \B\; \Rightarrow\; |\bar{q}(z)|\leq\delta^{-s},
\end{equation}
whence also
\begin{equation}\label{eq5}
|z|=1\;\Rightarrow\; |\bar{q}(z)|\leq\delta^{-s}.
\end{equation}
Now, the polynomial $p(z)=\prod_{i=1}^s(z-\theta_i)$ on the boundary of $\B$ clearly satisfies
$$
|p(z)|\geq \left[{1\over \bar{\delta}}-1\right]^s=\left[{1-\delta\over 2\delta}\right]^s,
$$
which combines with (\ref{eq10}) to imply that the modulus of the holomorphic in $\B$ function
$$
\bar{r}(z)=\left[\prod_{i=1}^s(\delta z-\theta_i)\right]^{-1}
$$
is bounded with $\delta^{-s}\left[{1-\delta\over 2\delta}\right]^{-s}=\left[{2\over 1-\delta}\right]^s$ on the boundary of $\B$. It follows that the coefficients $r_j$ of the Taylor series of $\bar{r}$ satisfy
$$
|r_j|\leq  \left[{2\over 1-\delta}\right]^s\bar{\delta}^{j},\;\;j=0,1,2,...
$$
When setting
\be
q^\ell(z)=p(z)r^\ell(z),\quad r^\ell(z)=\sum_{j=1}^\ell r_j z^j,
\ee{eq:qldiv}
for $|z|\leq 1$, utilizing the trivial upper bound $|p(z)|\leq 2^s$, we get
\be
|q^\ell(z)-\bar{q}(z)|
\leq |p(z)||r^\ell(z)-\bar{r}(z)|\leq 2^s\left[{2\over 1-\delta}\right]^s \sum_{j=\ell+1}^\infty |r_j|\leq \left[{4\over 1-\delta}\right]^s{\bar{\delta}^{\ell+1}\over 1-\bar{\delta}}.
\ee{eq12}
Note that $q^\ell(0)=p(0)r^\ell(0)=p(0)\bar{r}(0)=1$, that $q^\ell$ is a polynomial of degree $\ell+s$, and that $q^\ell$ is divisible by $p(z)$.
Besides this, on the unit circumference we have, by (\ref{eq12}),
\be
|q^\ell(z)|
\leq |\bar{q}(z)|+\left[{4\over 1-\delta}\right]^s{\bar{\delta}^{\ell+1}\over 1-\bar{\delta}}\leq \delta^{-s}+\underbrace{\left[{4\over 1-\delta}\right]^d{\bar{\delta}^{\ell+1}\over 1-\bar{\delta}}}_{R},
\ee{eq13}
where we used \eqref{eq5}.
Now,
$$
\bar{\delta}={2\delta\over1+\delta}={2-2\epsilon\over 2-\epsilon}={1-\epsilon\over 1-\epsilon/2}\leq 1-\epsilon/2\leq \e^{-\epsilon/2},
%=\exp\{-{\alpha\over 4ns}\}
$$
and
$$
{1\over 1-\bar{\delta}}={1+\delta\over 1-\delta}={2-\epsilon\over \epsilon}\leq {2\over\epsilon}.
$$
We can upper-bound $R$:
\[
R=\left[{4\over 1-\delta}\right]^s{\bar{\delta}^{\ell+1}\over 1-\bar{\delta}}
\leq {2^{2s+1}\over\epsilon^{s+1}}\e^{-\epsilon \ell/2}
%=\left[{4\over\epsilon}\right]^s
%\exp\{-{\alpha\over 4s}\}{2\over\epsilon}
%\leq \left[{4\over\epsilon}\right]^{s+1}\exp\{-{\alpha\over 4s}\}=\left[{8ns\over\alpha}\right]^{s+1}\exp\{-{\alpha\over 4s}\}.
\]Now, given positive integer $\ell$ and positive $\alpha$ such that
\begin{equation}\label{alphasmall}
{\alpha\over \ell}\leq \frac{1}{4},
\end{equation}
let $\epsilon={\alpha\over 2\ell s}$.
Since $0<\epsilon\leq {1\over 8}$, we have $-\log(\delta)=-\log(1-\epsilon)\leq 2\epsilon={\alpha\over \ell s},$ implying that $\bar{\delta}\leq \e^{-\epsilon/2}=\e^{-{\alpha\over 4\ell s}}$, and
\[
R\leq \left[{8\ell s\over\alpha}\right]^{s+1}\exp\left\{-{\alpha\over 4s}\right\}.
\]
Now let us put
\[%begin{equation}\label{alpha}
\alpha=\alpha(\ell,s)=4s(s+2)\log(8\ell s);
\]%end{equation}
observe that this choice of $\alpha$ satisfies (\ref{alphasmall}), provided that
\[%begin{equation}\label{nislarge}
\ell\geq O(1)s^2\log(s+1)
\]%end{equation}
with properly selected absolute constant $O(1)$. With this selection of $\alpha$, we have $\alpha\geq1$, whence
\begin{align*}
R\left[{\alpha\over \ell}\right]^{-1}
&\leq \exp\left\{-{\alpha\over 4s}\right\}\left[{8\ell s\over\alpha}\right]^{s+1}{\ell\over \alpha}\leq\exp\left\{-{\alpha\over 4s}\right\}[8\ell s]^{s+2}\nn
&\leq \exp\{-(s+2)\log(8\ell s)\}\exp\{(s+2)\log(8\ell s)\}=1,
\end{align*}
that is,
\begin{equation}\label{RR}
R\leq {\alpha\over \ell}\leq \frac{1}{4}.
\end{equation}
Furthermore,
\begin{equation}\label{firstterm}
\begin{array}{rcll}
\delta^{-s}&=&\exp\{-s\log(1-\epsilon)\}\leq \exp\{2\epsilon s\}=\exp\{{\alpha\over \ell}\}\leq 2,\\
\delta^{-2s}&=&\exp\{-2s\log(1-\epsilon)\}\leq \exp\{4\epsilon s\}=\exp\{{2\alpha\over \ell }\}\leq 1+\exp\{{1\over 2}\}{2\alpha\over \ell}\leq 1+{4\alpha\over \ell}.
\end{array}
\end{equation}
When invoking (\ref{eq13}) and utilizing (\ref{firstterm}) and (\ref{RR}) we get
\begin{align*}
{1\over 2\pi}\oint_{|z|=1}|q^\ell(z)|^2|dz|
&\leq \delta^{-2s}+2\delta^{-s}R+R^2\leq 1+4{\alpha\over \ell}+4R+{1\over 4}R \leq 1+10{\alpha\over \ell}.
\end{align*}
On the other hand, denoting by $q_0$, $q_1$,...,$q_{\ell+s}$ the coefficients of the polynomial $q^\ell$ and taking into account that $\bar{q}_0=q^\ell(0)=1$, we have
\begin{align}
\label{eq:qnorm}
1+\sum_{i=1}^{\ell+s}|q_i|^2=|q_0|^2+...+|q_{\ell+s}|^2
&={1\over 2\pi}\oint_{|z|=1}|q^\ell(z)|^2|dz|
\leq 1+10{\alpha\over \ell}.
%\eqno{\square}
\end{align}
We are done: when denoting $n=\ell+s$, and $q(z)=\sum_{i=1}^n q_j z^j$, we have the vector of coefficients $q=[0;q_1;...;q_n]\in \C^{n+1}$ of $q(z)$ such that, by \eqref{eq:qnorm},
\[
\|q\|^2_2\leq {40s(s+2)\log[8s(n-s)]\over n-s},
\]
and such that the polynomial $q^\ell(z)=1+q(z)$ is divisible by $p(z)$ due to \eqref{eq:qldiv}.
\qed

\subsection{Proof of Lemma~\ref{lem:sepsines}}
Let $\Pi_{\mathcal{S}_{2m}}$ be the $(2m+1)\times (2m+1)$ projector matrix built in the proof of Proposition \ref{prop:shift-inv-bilateral}, but now let $\phi^o \in \C_{m}(\Z)$  be obtained from the \textit{last} column of~$\Pi_{\mathcal{S}_{2m}}$.
As in that proof,  due to the shift-invariance of $\cH_s[\omega]$ we have   $x = \phi^o * \Delta^{m}x$ $\forall x\in \cH_s[\omega]$.  To prove the proposition
it remains to bound $\|\phi^o\|_2$.
\par
Note that in the premise of the proposition $\S_m$ is spanned by vectors
\begin{equation*}
\left\{ v(\omega):  [v(\omega)]_t = \frac{e^{i \omega_k t}}{\sqrt{m+1}}, \quad 0\leq t\leq 2m\right\}, \quad \omega  \in \{\omega_1, ..., \omega_s\}.
\end{equation*}
%see also Section~\ref{sec:intro-harmonic}.
Hence, the projector $\Pi_{\mathcal{S}_{2m}}$ can be written as
\[
\Pi_{\mathcal{S}_m} = V\left(V^\H V\right)^{-1}V^\H,
\]
where $V$ is an $(2m+1)\times s$ Vandermonde matrix  with columns $v(\omega_k)$, $k= 1,...,s$.
Note that since $s \le 2m+1$, and $\omega_k$, $k = 1,...,s$ are distinct, matrix~$V$ has full column rank.
Now, in order to bound $\|\phi^o\|_2$ from above it suffices to separate the minimal eigenvalue  $\lmin(V^\H V)$ of  $V^\H V$ from zero.
Indeed, assuming that $\lmin(V^\H V)  > 0$ we may write
\begin{align*}
\Pi_{\mathcal{S}_m} = U U^\H,
\end{align*}
where $U = [U_1, ..., U_s]$ is the unitary normalization of $V$:
\[
U = [U_1 \cdots U_s] = V(V^\H V)^{-1/2}, \quad U^\H U = I_s.
\]
Let $u = [u_1, ..., u_s]$ be the last row of $U$, and $v$ that of $V$. Note that the vector $\psi=uU^\H = \sum_{k=1}^s u_k [U_k]^\H$ has the same $\ell_2$-norm as $\phi^o$, and so
 $\|\phi^o\|_2^2 = \|u\|_2^2$.
On the other hand, because
$
u =v (V^\H V)^{-1/2},
$
we arrive at
\begin{align*}
\|u\|_{2}^2 \le \|v\|_2^2\lmin^{-1}(V^\H V) \le \frac{s}{2m+1}\lmin^{-1}(V^\H V)
\end{align*}
where the last inequality is due to the bound $(2m+1)^{-1/2}$ on the moduli of elements of $v$.
Finally, we utilize the bound on the condition number of a Vandermonde matrix:% (see~\cite{moitra2015super}):
\begin{lemma}[{\cite[Theorem~2.3]{moitra2015super}}]
Let $\delta_{\min}$ be given by~\eqref{def:min-sep}; one has
\[
\frac{\lmax(V^\H V)}{\lmin(V^\H V)} \le \left({m - \frac{2\pi}{\delta_{\min}}}\right)^{-1} \left({m + \frac{2\pi}{\delta_{\min}}}\right).
\]
\end{lemma}
\noindent We clearly have $\| V \|_* \ge 1$, whence $\lmax(V^\H V) \ge 1$. Together with~\eqref{eq:min-sep} this results in
\[
\lmin^{-1}(V^\H V) \le \frac{\nu + 1}{\nu- 1},
\]
whence the required bound on $\|\phi^o\|_2$.
\qed
\subsection{Proof of Proposition \ref{pr:sines}}
Note that in the premise of the proposition $k=\lfloor L/(s\log[L])\rfloor$ is correctly defined and $K=L-2k\geq L/2$ so that
\be
\kappa_{K,k}\leq C(s\log L)^{-1/2} \quad\mbox{and}\quad \kappa_{k,K}\leq C'\sqrt{s\log L}.
\ee{eq:mykap}
%\paragraph{1$^o$.}
When applying Proposition~\ref{pr:l2sloss} (recall that $\varkappa=0$ in our setting),
we conclude that the error of the estimate $\wh \vphi *y$ satisfies, with probability at least $1-\alpha/3$,
\be
\|x-\wh x \|_{K,2}
\le C \sigma \left(\kappa_{k,K}s\sqrt{\log[1/\alpha]} + \sqrt{s\log [L/\alpha]}\right).
\ee{eq:sines-best-center}
%
%\paragraph{2$^o$.}
On the other hand,  due to $\kappa_{K,k}\leq 1$, applying  Proposition \ref{pr:l2sloss-p} we conclude that with probability $1-\alpha/3$  the error of the left estimate
$\wh\vphi^+*\Delta^my$  satisfies:
\[\left\|\Delta^{-m}(x - \wh\vphi^+ * \Delta^my)\right\|_{k,2}
\le C'\sigma\left(\kappa_{K,k} s^2\log[L]\sqrt{\log[1/\alpha]}+s\sqrt{\log[L]\log[L/\alpha]}\right),
\]
and the same estimation holds true for the right estimate $\wh\vphi^-*\Delta^{-m}y$:
\[\left\|\Delta^{m}(x - \wh\vphi^- * \Delta^{-m}y)\right\|_{k,2}
\le C'\sigma\left(\kappa_{K,k} s^2\log[L]\sqrt{\log[1/\alpha]}+s\sqrt{\log[L]\log[L/\alpha]}\right).
\]

When combining the latter bounds with \rf{eq:sines-best-center} we arrive at the bound with probability $\geq 1-\alpha$:
\bse
\left\|x - \wh x\right\|_{L,2}&\leq &\left\|\Delta^{m}(x - \wh\vphi^- * \Delta^{-m}y)\right\|_{k,2}+\|x-\wh x \|_{K,2}+\left\|\Delta^{-m}(x - \wh\vphi^+ * \Delta^my)\right\|_{k,2}\\
&\leq & %C \sigma \left(\kappa_{k,K}s\sqrt{\log[1/\alpha]} + \sqrt{s\log [L/\alpha]}\right)+2 C'\sigma\left(\kappa_{K,k} s^2\log[L]\sqrt{\log[1/\alpha]}+s\sqrt{\log[L]\log[L/\alpha]}\right)
C\sigma s\sqrt{\log[L]\log[L/\alpha]}+C'\sigma s\sqrt{\log[1/\alpha]}(\kappa_{k,K}+\kappa_{K,k}s\log [L])\\
\mbox{(by $\rf{eq:mykap}$)}&\leq &C\sigma s\log[L/\alpha]+C''\sigma s\sqrt{s\log[L]\log[1/\alpha]}\leq C\sigma s^{3/2}\log[L/\alpha].\qquad\qquad\mbox{\qed}
\ese
\section{Naive adaptive estimate}\label{sec:naive}
In this section,\footnote{We use notation defined in Sections \ref{sec:intro-notation} and \ref{sec:addnot}.} we consider the ``naive'' adaptive estimate $\wh x=\wh \phi* y$ where $\wh \phi\in \C_m(\Z)$ solves the optimization problem
\be
\min_{\phi\in \C_m(\Z)}\|y-\phi*y\|_{n,2} \; \mbox{ subject to } \; \|\phi\|_2\leq {\rho\over \sqrt{2m+1}}.
\ee{naive1}
Recall that our goal is to show that using estimate $\wh x$ is really not  a good idea.  To make the long story short, from now on, we consider the simplified version of the estimation problem in which $m=n$, signals are $2m+1$-periodic, % on $-m,...,m$
and linear estimates are in the form of circular (periodic) convolution
\[
[\phi*y]_t=\sum_{\tau=-m}^m\phi_{\tau}y_{s(t,\tau)},\qquad |t|\leq m,
\]
where $s(t,\tau)=[t+m-\tau \ \mathrm{mod}\ 2m+1]-m$.
Because the Discrete Fourier Transform diagonalizes the periodic convolution, problem \rf{naive1} may be equivalently reformulated in the space of Fourier coefficients
\be
\min_{w\in {\C^{2m+1}}}\|z-Zw\|_{n,2} \; \mbox{ subject to } \; \|w\|_2\leq {\rho}
\ee{naivef}
where $z=F_m[y]$, $Z=\diag(z)$ (with $A = \diag(a)$ being the diagonal matrix with entries $A_{ii}=a_i$), and $w$ is a properly ``rephased'' DFT of $\phi$ with $|w_k|=\sqrt{2m+1}\,|(F_m[\phi])_k|$, $1\leq k\leq 2m+1$.
\par
Consider the situation in which the signal to recover is just one ``complex sinusoid,'' e.g., $x_\tau=ae^{2\pi i\tau\over 2m+1}$, $\tau \in \Z$, $a\in \C$, and let us show that the error of the naive estimate may be much larger than the ``oracle'' error. We have $F_m[x]=fe_1$ where $e_1$ is the first basis orth,  $f=a\sqrt{2m+1}$ with $|f|=\|x\|_{m,2}=|a|\sqrt{2m+1}$, and the ``sequence-space'' observation $z$ satisfies
\[
z=fe_1+\sigma \zeta, \quad \zeta\sim \CN(0,I_n).
\]
Obviously, in this case there exist a filter $\phi^o$ with $\|\phi^o\|_2=(2m+1)^{-1/2}$ such that $x=\phi^o*x$, so that the integral $\alpha$-risk of the ``oracle estimate'' $\phi^o*y$ is $\cO(\sigma)$ up to logarithmic in $\alpha$ factor. Let us show that in this simple situation the risk of the naive estimate may be significantly higher.
\par
First of all, note that the optimal solution $\wh w$ to the problem \rf{naivef} with $\rho=1$ is of the form
\[
\wh w_k={|z_k|^2\over |z_k|^2+\lambda},\;\;1\leq k\leq 2m+1
\]
where $\lambda$ is chosen to ensure $\|\wh w\|_2=1$. Let us bound $\lambda$ from below. We have
\bse
1=\|\wh w\|_2^2={|z_1|^4\over (|z_1|^2+\lambda)^2}+\sum_{k=2}^{2m+1}{\sigma^4 |\zeta_k|^4\over (\sigma^2|\zeta_k|^2+\lambda)^2}
\geq \sum_{k=2}^{2m+1}{\sigma^4 |\zeta_k|^4\over (\sigma^2M_m^2+\lambda)^2}\geq {\sigma^4S_m^2\over 2m(\lambda+\sigma^2M_m)^2}
\ese
where $M_m=\max_{1\leq k\leq 2m+1}|\zeta_k|^2$ and $S_m=\sum_{k=2}^{2m+1}|\zeta_k|^2$.
Since  with {high probability} (say, $1-\cO(1/m)$) $M_m=\cO(\log m)$ and $S_m=\cO(m)$ (cf.~\rf{eq:complex_gaussian_max} and {\rf{eq:chi_bound_low}}), for $m$ large enough one has
\[
\lambda\geq \sigma^2\left({S_m\over \sqrt{2m}}-M_m\right)\geq c\sigma^2 \sqrt{m}
\]
with probability at least $1-\cO(1/m)$.
As a result,
\[
1-\wh w_1=1-{|z_1|^2\over |z_1|^2 +\lambda}={\lambda\over |z_1|^2 +\lambda}\geq {\lambda\over (|f|+\sigma|M_n|)^2+\lambda}\geq c'
\]
whenever $f$ satisfies $|f|^2\leq C\sigma^2 \sqrt{m}$. Next, observe that
\bse
\|x-\wh x\|^2_{m,2}&=&
\|F_m[x]-Z\wh w\|_2^2=\|fe_1-Z\wh w\|_2^2\geq |f-z_1\wh w_1|^2\\&\geq& \half |f(1-\wh w_1)|^2-\sigma^2|\zeta_1|^2\wh w_1^2\geq
c|f|^2- \sigma^2M_m\geq c'|f|^2
\ese
for $|f|\geq c''\sigma \sqrt{\log m}$. In other words, when the signal amplitude satisfies
\[
{c\sigma^2\log m \over m}\leq |a|^2\leq {C\sigma^2 \over \sqrt{m}},
\]
the loss $\|\wh x-x\|_{m,2}$ of the naive estimate is lower bounded, with probability at least $1-\cO(1/m)$, with $c'\|x\|_{m,2}$. In particular, when $a\asymp{\sigma m^{-1/4}}$ this error is at least order of $\sigma m^{1/4}$, which is incomparably worse than the error $\cO(\sigma)$ of the oracle estimate.

\bibliography{references}
\bibliographystyle{plain}

\end{document}

%% file: preamble.tex
\usepackage[colorlinks]{hyperref}
\usepackage[dvipsnames]{xcolor}
\hypersetup{
    linkcolor=red,
    citecolor=OliveGreen,
    filecolor=black,
    urlcolor=black,
}
\usepackage{caption}

\newcommand{\Win}{w}
\newcommand{\beq}[1]{\begin{equation}\label{eq:#1}}
\newcommand{\eeq}{\end{equation}}
\newcommand{\eref}[1]{\eqref{eq:#1}}
\newcommand{\con}{\textup{con}}
\newcommand{\pen}{\textup{pen}}
\newcommand{\Rem}{\textup{Q}}

\renewcommand{\Re}{\operatorname{Re}}
\renewcommand{\Im}{\operatorname{Im}}
\usepackage{verbatim}
\usepackage{amsthm,amsmath}
\usepackage{amssymb}
\usepackage{bm}
\usepackage{bbm}
\usepackage{dsfont}
\usepackage{comment}
\usepackage{multirow}
\usepackage{color}
\usepackage{tikz}
\usepackage{mathtools}
\usepackage{pifont}
\usepackage{amsmath}
\usepackage{enumerate, algorithm,algorithmic,caption,graphicx,color}
\usepackage{chngcntr}

\newcommand{\half}{ \mbox{\small$\frac{1}{2}$}}

\newcommand{\be}{\begin{eqnarray}}
\newcommand{\ee}[1]{\label{#1}\end{eqnarray}}
\newcommand{\nn}{\nonumber \\}
\newcommand{\ese}{\end{eqnarray*}}
\newcommand{\bse}{\begin{eqnarray*}}
\newcommand{\rf}[1]{\eqref{#1}}
\def\qed{\hfill$\square$}

\def\cX{{\cal X}}
\def\cS{{\cal S}}

\def\cO{O}

\def\cH{{\mathcal{H}}}
\def\X{{\mathcal{X}}}

\def\pen{{\hbox{\scriptsize\rm pen}}}

\def\Tr{\mathop{\hbox{\rm Tr}}}

\def\Diag{{\hbox{\rm Diag}}}

\def\E{{\mathbf{E}}}

\def\N{{\mathcal{N}}}
\def\CN{{\C\mathcal{N}}}
\def\S{{\mathcal{S}}}
\def\F{{\mathcal{F}}}
\def\C{\mathds{C}}
\def\R{\mathds{R}}
\def\Z{\mathds{Z}}

\newcommand{\Prob}{\mathrm{Prob}}
\def\e{{\rm e}}
\newcommand{\wh}[1]{{\widehat{#1}}}

\def\lmax{\lambda_{\max}}
\def\lmin{\lambda_{\min}}

\def\F{{\textup{F}}}
\def\T{{\textup{T}}}
\def\H{{\textup{H}}}
\def\vphi{\varphi}
\def\vrho{\varrho}

\newcommand{\prob}{\mathrm{Prob}}
\def\risk{{\hbox{\rm Risk}}}

\def\DirKer{{\cal D}}

\newtheorem{theorem}{Theorem}[section]
\newtheorem{definition}{Definition}[section]
\newtheorem{ass}{Assumption}[section]

\newtheorem{lemma}{Lemma}[section]

\newtheorem{prop}{Proposition}[section]

\allowdisplaybreaks